\documentclass[reqno,11pt]{article}
\usepackage{mathrsfs}
\usepackage{amssymb}
\usepackage{amsthm}
\usepackage{amsmath}
\usepackage{amsfonts}
\usepackage{enumerate}
\usepackage{bm, enumerate,xcolor}

\usepackage{graphicx}

\numberwithin{equation}{section}

\newtheorem{theorem}{Theorem}[section]
\newtheorem{lemma}[theorem]{Lemma}
\newtheorem{corollary}[theorem]{Corollary}
\newtheorem{proposition}[theorem]{Proposition}

\theoremstyle{definition}
\newtheorem{definition}[theorem]{Definition}
\newtheorem{assumption}[theorem]{Assumption}
\newtheorem{example}[theorem]{Example}

\theoremstyle{remark}
\newtheorem{remark}[theorem]{Remark}

\begin{document}
\title{On Some Smooth Symmetric Transonic Flows with Nonzero Angular Velocity and Vorticity}
\author{Shangkun WENG\thanks{School of mathematics and statistics, Wuhan University, Wuhan, Hubei Province, 430072, People's Republic of China. Email: skweng@whu.edu.cn}\and Zhouping XIN\thanks{The Institute of Mathematical Sciences, The Chinese University of Hong Kong, Shatin, NT, Hong Kong. E-mail: zpxin@ims.cuhk.edu.hk}\and Hongwei YUAN\thanks{The Institute of Mathematical Sciences, The Chinese University of Hong Kong, Shatin, NT, Hong Kong. E-mail: hwyuan@link.cuhk.edu.hk}}
\date{}
\maketitle

\def\be{\begin{eqnarray}}
\def\ee{\end{eqnarray}}
\def\ba{\begin{align}}
\def\ea{\end{align}}
\def\bay{\begin{array}}
\def\eay{\end{array}}
\def\bca{\begin{cases}}
\def\eca{\end{cases}}
\def\p{\partial}
\def\hphi{\hat{\phi}}
\def\bphi{\bar{\phi}}
\def\tphi{\tilde{\phi}}
\def\no{\nonumber}
\def\eps{\epsilon}
\def\veps{\varepsilon}
\def\de{\delta}
\def\De{\Delta}
\def\om{\omega}
\def\Om{\Omega}
\def\f{\frac}
\def\th{\theta}
\def\vth{\vartheta}
\def\la{\lambda}
\def\lab{\label}
\def\b{\bigg}
\def\var{\varphi}
\def\na{\nabla}
\def\ka{\kappa}
\def\al{\alpha}
\def\La{\Lambda}
\def\ga{\gamma}
\def\Ga{\Gamma}
\def\ti{\tilde}
\def\wti{\widetilde}
\def\wh{\widehat}
\def\ol{\overline}
\def\ul{\underline}
\def\Th{\Theta}
\def\si{\sigma}
\def\Si{\Sigma}
\def\oo{\infty}
\def\q{\quad}
\def\z{\zeta}
\def\co{\coloneqq}
\def\eqq{\eqqcolon}
\def\di{\displaystyle}
\def\bt{\begin{theorem}}
\def\et{\end{theorem}}
\def\bc{\begin{corollary}}
\def\ec{\end{corollary}}
\def\bl{\begin{lemma}}
\def\el{\end{lemma}}
\def\bp{\begin{proposition}}
\def\ep{\end{proposition}}
\def\br{\begin{remark}}
\def\er{\end{remark}}
\def\bd{\begin{definition}}
\def\ed{\end{definition}}
\def\bpf{\begin{proof}}
\def\epf{\end{proof}}
\def\bex{\begin{example}}
\def\eex{\end{example}}
\def\bq{\begin{question}}
\def\eq{\end{question}}
\def\bas{\begin{assumption}}
\def\eas{\end{assumption}}
\def\ber{\begin{exercise}}
\def\eer{\end{exercise}}
\def\mb{\mathbb}
\def\mbR{\mb{R}}
\def\mbZ{\mb{Z}}
\def\mc{\mathcal}
\def\mcS{\mc{S}}
\def\ms{\mathscr}
\def\lan{\langle}
\def\ran{\rangle}
\def\lb{\llbracket}
\def\rb{\rrbracket}
\def\fr#1#2{{\frac{#1}{#2}}}
\def\dfr#1#2{{\dfrac{#1}{#2}}}
\def\u{{\textbf u}}
\def\v{{\textbf v}}
\def\w{{\textbf w}}
\def\d{{\textbf d}}
\def\nn{{\textbf n}}
\def\x{{\textbf x}}
\def\e{{\textbf e}}
\def\D{{\textbf D}}
\def\U{{\textbf U}}
\def\M{{\textbf M}}
\def\F{{\mathcal F}}
\def\I{{\mathcal I}}
\def\W{{\mathcal W}}
\def\div{{\rm div\,}}
\def\curl{{\rm curl\,}}
\def\R{{\mathbb R}}
\def\FF{{\textbf F}}
\def\A{{\textbf A}}
\def\R{{\textbf R}}
\def\r{{\textbf r}}

\begin{abstract}

  This paper concerns the structural stability of smooth cylindrically symmetric transonic flows in a concentric cylinder. Both cylindrical and axi-symmetric perturbations are considered. The governing system here is of mixed elliptic-hyperbolic and changes type and the suitable formulation of boundary conditions at the boundaries is of great importance. First, we establish the existence and uniqueness of smooth cylindrical transonic spiral solutions with nonzero angular velocity and vorticity which are close to the background transonic flow with small perturbations of the Bernoulli's function and the entropy at the outer cylinder and the flow angles at both the inner and outer cylinders independent of the symmetric axis, and it is shown that in this case, the sonic points of the flow are nonexceptional and noncharacteristically degenerate, and form a cylindrical surface. Second, we also prove the existence and uniqueness of axi-symmetric smooth transonic rotational flows which are adjacent to the background transonic flow, whose sonic points form an axi-symmetric surface. The key elements in our analysis are to utilize the deformation-curl decomposition for the steady Euler system introduced in \cite{WengXin19} to deal with the hyperbolicity in subsonic regions and to find an appropriate multiplier for the linearized second order mixed type equations which are crucial to identify the suitable boundary conditions and to yield the important basic energy estimates.

\end{abstract}

\begin{center}
\begin{minipage}{5.5in}
Mathematics Subject Classifications 2010: 76H05, 35M12, 35L65, 76N15.\\
Key words: smooth transonic spiral flows, mixed type, multiplier, deformation-curl decomposition, sonic surface.
\end{minipage}
\end{center}
\section{Introduction and main results}\noindent

In this paper, we study the structural stability of some smooth cylindrically symmetric transonic spiral flows in a concentric cylinder $\tilde{\Omega}=\{(x_1,x_2,x_3): r_0<r=\sqrt{x_1^2+x_2^2}<r_1,x_3\in \mathbb{R}\}$. The flow is governed by the following steady compressible Euler system
\begin{align}\label{comeuler3d}
\begin{cases}
\partial_{x_1}(\rho u_1)+\partial_{x_2}(\rho u_2)+\partial_{x_3}(\rho u_3)=0,\\
\partial_{x_1}(\rho u_1^2)+\partial_{x_2}(\rho u_1 u_2)+\partial_{x_3}(\rho u_1 u_3)+\partial_{x_1} p=0,\\
\partial_{x_1}(\rho u_1u_2)+\partial_{x_2}(\rho u_2^2)+\partial_{x_3}(\rho u_2 u_3)+\partial_{x_2} p=0,\\
\partial_{x_1}(\rho u_1u_3)+\partial_{x_2}(\rho u_2 u_3)+\partial_{x_3}(\rho u_3^2)+\partial_{x_3} p=0,\\
\displaystyle\sum_{j=1}^3\partial_{x_j}(\rho u_j (e+\frac{1}{2}|{\bf u}|^2)+u_j p)=0,\\
\end{cases}
\end{align}
where ${\bf u}=(u_1,u_2, u_3)^t$, $\rho, p$ are the velocity, the density and pressure respectively, $e=e(\rho,p)$ represents the internal energy. Here we consider only the polytropic gas, therefore $p= A(S) \rho^{\gamma}$, where $A(S)= a e^{S}$ and $\gamma\in (1,3), a$ are positive constants, $e=\frac{p}{(\gamma-1)\rho}$. Denote the Bernoulli's function and the local sonic speed by $B=\frac12|{\bf u}|^2+\frac{\gamma p}{(\gamma-1)\rho}$ and $c(\rho, A)=\sqrt{A\gamma} \rho^{\frac{\gamma-1}{2}}$, respectively.

Introduce the cylindrical coordinates $(r, \theta, x_3)$
\begin{eqnarray}\no
r=\sqrt{x_1^2+x_2^2},\ \theta=\arctan \frac{x_2}{x_1},\ \ x_3=x_3
\end{eqnarray}
and decompose the velocity as ${\bf u}= U_1 {\bf e}_r + U_2 {\bf e}_{\theta} +U_3 {\bf e}_3$ with
\begin{align}\no
{\bf e}_r=\begin{pmatrix} \cos\theta\\ \sin\theta\\ 0\end{pmatrix},
{\bf e}_{\theta}=\begin{pmatrix}-\sin\theta\\ \cos\theta\\ 0\end{pmatrix},  {\bf e}_{3}=\begin{pmatrix} 0\\ 0\\ 1\end{pmatrix}.
\end{align}
Then the system \eqref{comeuler3d} can be rewritten as
\begin{eqnarray}\label{ProblemImm}
\begin{cases}
\partial_r(\rho U_1)+\frac{1}{r}\partial_{\theta}(\rho U_2)+\frac{1}{r} \rho U_1+ \partial_3(\rho U_3)=0,\\
(U_1\partial_r +\frac{U_2}{r}\partial_{\theta}+U_3 \partial_3) U_1+\frac{1}{\rho} \partial_r p-\frac{U_2^2}r=0,\\
(U_1\partial_r +\frac{U_2}{r}\partial_{\theta}+U_3 \partial_3) U_2+\frac{1}{r\rho} \partial_{\theta}p+\frac{U_1U_2}r=0,\\
(U_1\partial_r +\frac{U_2}{r}\partial_{\theta}+U_3 \partial_3) U_3+\frac{1}{\rho} \partial_{3}p=0,\\
(U_1\partial_r +\frac{U_2}{r}\partial_{\theta}+U_3 \partial_3) A=0.
\end{cases}
\end{eqnarray}

In the cylindrical coordinates, the vorticity has the form $\text{curl }{\bf u}=\omega_r {\bf e}_r + \omega_{\theta}{\bf e}_{\theta}+ \omega_3 {\bf e}_3$, where
\begin{eqnarray}\label{vor1}
\omega_r=\frac{1}{r}\partial_{\theta} U_3 -\partial_3 U_2,\ \ \omega_{\theta}= \partial_3 U_1- \partial_r U_3,\ \omega_3=\frac{1}{r}\partial_r(r U_2)-\frac{1}{r}\partial_{\theta} U_1.
\end{eqnarray}

Courant and Friedrichs in \cite[Section 104]{Courant1948} had used the hodograph method to obtain some special planar radially symmetric flows including circulatory flows, purely radial flows and their superimposing (spiral flows). Circulatory flows and purely radial flows are radially symmetric flows with only angular and radial velocity, respectively. These radial symmetrical flows can be regarded as cylindrical symmetric flows in three-dimensional setting. It had been proved in \cite{Courant1948} that spiral flows can take place only outside a limiting circular cylinder where the Jacobian of the hodograph transformation is zero and may change smoothly from subsonic to supersonic or vice verse. We study the cylindrical symmetric smooth solutions in \cite{WXY20a} by considering one side boundary value problem to the steady Euler system in a cylinder and analyze the dependence of the solutions on the boundary data, which suits our purpose for the investigation of structural stability of this special transonic flows in this paper. Here we deal with structural stability of a special class of smooth cylindrically symmetric, irrotational transonic spiral flows moving from the outer into the inner cylinders. More precisely, the background flow is described by smooth functions of the form ${\bf u}(x)= U_{1}(r) {\bf e}_r + U_{2}(r) {\bf e}_{\theta}, \rho(x)=\rho(r)$ and $p(x)=p(r), A(x)=A(r)$, solving the following system
\begin{align}\label{radial-euler}
\begin{cases}
\frac{d}{dr}(\rho U_1)+\frac1r \rho U_1=0, \ \ &0<r_0<r<r_1,\\
U_1U_1'+\frac1\rho \frac{d}{dr}p-\frac{U_2^2}r=0,  \ \ &0<r_0<r<r_1,\\
U_1U_2'+\frac{U_1U_2}r=0, \ \ &0<r_0<r<r_1,\\
U_1A'=0,  \ \ &0<r_0<r<r_1
\end{cases}
\end{align}
with the boundary conditions at the outer cylinder $r=r_1$:
\begin{eqnarray}\label{ProblemI}
\rho(r_1)=\rho_{0}>0,\ U_{1}(r_1)=U_{10}\leq 0, \ U_{2}(r_1)=U_{20}\neq 0,\ A(r_1)=A_{0}>0.
\end{eqnarray}
Let $B_0=\frac{1}{2}(U_{10}^2+ U_{20}^2)+ \frac{\gamma}{\gamma-1} A_0 \rho_0^{\gamma-1}$ and denote the Mach numbers by
\begin{eqnarray}\nonumber
M_{b1}(r)=\frac{U_{b1}(r)}{c(\rho_b, A_0)},\ \ M_{b2}(r)=\frac{U_{b2}(r)}{c(\rho_b, A_0)},\ \ {\bf M}_b(r)= (M_{b1}(r), M_{b2}(r))^t.
\end{eqnarray}


The following proposition had been established in \cite{WXY20a}.

\begin{proposition}\label{background}
{\it Suppose that the incoming flow is subsonic, i.e. $A_0\gamma \rho_0^{\gamma-1}> U_{10}^2+U_{20}^2$. Then there exist constants $0<r^\sharp<r_c<r_1$, where $r^\sharp$ depends only on $r_1$, $\ga$ and the incoming flow at $r_1$,
\begin{eqnarray}\nonumber
r_c=\sqrt{\frac{(\gamma+1)(\kappa_1^2+ \kappa_2^2\rho_c^2)}{2(\gamma-1)B_0\rho_c^2}},\ \ \rho_c=\bigg(\frac{2(\ga-1)B_0}{(\ga+1)\ga A_0}\bigg)^{\frac1{\ga-1}},\ \  \kappa_1=r_1\rho_0U_{10}, \ \ \kappa_2=r_1U_{20},
\end{eqnarray}
such that if $r^\sharp<r_0<r_c$, there exists a unique smooth irrotational transonic spiral flow $(U_{b1}, U_{b2}, \rho_b, A_0)$ to \eqref{radial-euler}-\eqref{ProblemI} in $[r_0,r_1]$ with all sonic points located at the cylinder $r=r_c$, which moves from the outer cylinder to the inner one. The total Mach number $|{\bf M}_b(r)|$ monotonically increases as $r$ decreases from $r_1$ to $r_0$, and $|{\bf M}_b(r)|<1$ for any $r\in (r_c,r_1]$ and $|{\bf M}_b(r)|>1$ for any $r\in [r_0,r_c)$, but the radial Mach number $|M_1(r)|$ is always less than $1$ for any $r\in [r_0,r_1]$. Moreover, all the sonic points are nonexceptional and noncharacteristically degenerate.
}\end{proposition}



In particular, if $U_{10}=0$, then $U_1\equiv 0, A(r)\equiv A_0$, $U_2(r)=\frac{\kappa_2}{r}$ and
\begin{eqnarray}\nonumber
\rho(r)=\left(\frac{\gamma-1}{A_0\gamma}\right)^{\frac{1}{\gamma-1}} \left(B_0-\frac{\kappa_2^2}{2r^2}\right)^{\frac{1}{\gamma-1}},\ \ r^{\sharp}=\frac{|\kappa_2|}{\sqrt{2B_0}}.
\end{eqnarray}
All the sonic points locate at the cylinder $r=r_c=\sqrt{\frac{\gamma+1}{2(\gamma-1) B_0}}|\kappa_2|$. Also if $U_{10}<0$, then $U_1(r)<0$ for any $r\in [r_0,r_1]$. Fix $r_0\in(r^\sharp,r_c)$, we call this smooth transonic spiral flow constructed in Proposition \ref{background} on $[r_0,r_1]$ to be the background flow. Note that the background flow is always irrotational.

Our main concern in this paper is the structural stability of this background flow under suitable perturbations of the boundary data. Note that the existence and stability of subsonic circulatory flows outside a smooth profile have been studied extensively by \cite{bers54,Bers1958,Cui2011,fg57,shiffman52}. Yet due to the degeneracy, the existence and structural stability of smooth transonic flows are substantially difficult. Following the definition by Bers \cite{Bers1958}, a sonic point in a $C^2$ transonic flow is exceptional if and only if the velocity is orthogonal to the sonic curve at this point. Due to the nonzero angular velocity, all the sonic points of the background flow are nonexceptional and noncharacteristically degenerate. This is quite different from the recent important results obtained by Wang and Xin in \cite{WX2013,WX2019,WX2020}, where the existence and uniqueness of smooth transonic flow of Meyer type in de Laval nozzles are proved and all the sonic points on the throat are exceptional and characteristically degenerate. In this paper, we will examine the structural stability of this kind of background flows by prescribing some physically acceptable boundary conditions at the entrance and exit and establish the existence and uniqueness of two classes of smooth transonic spiral flows with small nonzero vorticity, whose sonic points are all nonexceptional and noncharacteristic.

It seems quite difficult to analyze the structural stability of the background flow under generic three dimensional perturbations. Therefore, we consider here two classes of perturbations with special symmetries. First we investigate the structural stability under a perturbation which is independent of $x_3$. Thus we will identify a class of physical acceptable boundary conditions at the circular cylinders respectively, and establish the existence and uniqueness of a class of smooth transonic spiral flows with small nonzero vorticity which are independent of the symmetry axis and close to the background flow. More precisely, we are looking for solutions to \eqref{ProblemImm} with the form $(U_1(r,\theta), U_2(r,\theta), U_3(r,\theta)\equiv 0, \rho(r,\theta), p(r,\theta))$ satisfying
\be\label{2d}\begin{cases}
\partial_r(\rho U_1)+\frac{1}{r}\partial_{\theta}(\rho U_2)+\frac{1}{r} \rho U_1=0,\\
(U_1\partial_r +\frac{U_2}{r}\partial_{\theta}) U_1+\frac{1}{\rho} \partial_r p-\frac{U_2^2}r=0,\\
(U_1\partial_r +\frac{U_2}{r}\partial_{\theta}) U_2+\frac{1}{r\rho} \partial_{\theta}p+\frac{U_1U_2}r=0,\\
(U_1\partial_r +\frac{U_2}{r}\partial_{\theta}) A=0
\end{cases}\ee
with
\be\label{2dboundary1}
&&B(r_1,\theta)=B_0 + \epsilon B_1(\theta),\\\label{2dboundary2}
&&A(r_1,\theta)=A_0 + \epsilon A_1(\theta),\\\label{2dboundary3}
&&U_1(r_0,\theta)-l_0 U_2(r_0,\theta)=a_0+\epsilon g_0(\theta),\\\label{2dboundary4}
&&U_2(r_1,\theta)=a_1+\epsilon g_1(\theta),
\ee
where $B_1, A_1, g_0$ and $g_1$ are given periodic functions with period $2\pi$. $l_0$ is a constant in a suitable range to be specified later and
\begin{eqnarray}\nonumber
a_0= U_{b1}(r_0)- l_0 U_{b2}(r_0),\ \ \ a_1=U_{20}.
\end{eqnarray}

Since the Bernoulli's quantity and the entropy satisfy transport equations, it is natural to prescribe the boundary conditions \eqref{2dboundary1} and \eqref{2dboundary2} at the entrance. We prescribe some restrictions on the flow angles \eqref{2dboundary3} and \eqref{2dboundary4} at the entrance and exit, which are physically acceptable and experimentally controllable. The flow angle restrictions \eqref{2dboundary3} and \eqref{2dboundary4} are also admissible for the linearized mixed type potential equation from the mathematical point of view (see Lemma \ref{basic-energy}), and are helpful to yield the important basic energy estimates.

Clearly, due to the independence of $x_3$, the problem can be reduced to the stability theory for the background flow in the annulus $\Omega=\{(r,\theta): r_0<r<r_1,\theta\in [0,2\pi]\}$.

{\bf Notation.} $\|\cdot\|_{k}$ will denote the norm in the Sobolev space $H^{k}(\Omega)$ for $k=1,2,3,4$. $\|\cdot\|_{L^p}$ will be the norm in $L^p(\Omega)$ for any $1\leq p\leq \infty$. Note that the estimate $\|uv\|_1\leq C\|u\|_2\|v\|_1$ is true if $u\in H^2(\Omega)$, $v\in H^1(\Omega)$. $\mathbb{T}_{2\pi}$ denotes the 1-d torus with period $2\pi$.

Then the first main result can be stated as follows.
\begin{theorem}\label{2dmain}
{\it Assume that $\gamma\in (1,3)$. Let a background flow with nonzero radial velocity $U_{b1}\neq 0$ be given such that
\begin{eqnarray}\label{incoming}
4-(3-\gamma)|{\bf M}_b(r_0)|^2>0.
\end{eqnarray}
Then for any constant $l_0$ with
\be\label{bo}
l_0\not\in \left(\frac{M_{b1}(r_0)M_{b2}(r_0)-\sqrt{|{\bf M}_{b}(r_0)|^2-1}}{1-M_{b1}^2(r_0)},\frac{M_{b1}(r_0)M_{b2}(r_0)+\sqrt{|{\bf M}_{b}(r_0)|^2-1}}{1-M_{b1}^2(r_0)}\right),
\ee
and any boundary data $g_0\in H^3(\mathbb{T}_{2\pi}), g_1\in C^{3,\alpha}(\mathbb{T}_{2\pi})$ and $(B_1, A_1)\in (C^{4,\alpha}(\mathbb{T}_{2\pi}))^2$ for some $\alpha\in (0,1)$, there exists a small constant $\epsilon_0$ depending on the background flow, $l_0$ and boundary datum $g_0, g_1, B_1, A_1$, such that for any $0<\epsilon<\epsilon_0$, the problem \eqref{2d} has a unique smooth transonic solution with possible nonzero vorticity $(U_1,U_2, B, A)\in (H^3(\Omega))^2\times(H^4(\Omega))^2$, which satisfies the boundary conditions \eqref{2dboundary1}-\eqref{2dboundary4} and the estimate
\be\label{2d10}
\|U_1-U_{b1}\|_3+\|U_2-U_{b2}\|_3+\|B-B_{0}\|_4+\|A-A_0\|_4\leq C\epsilon,
\ee
for some constant $C$ depending only on the background flow and the boundary datum.

Moreover, all the sonic points form a closed arc with a parametric representation $r=s(\theta)\in C^{1}(\mathbb{T}_{2\pi})$. The sonic curve is closed to the background sonic circle in the sense that
\be\label{2d2}
\|s(\theta)-r_c\|_{C^1(\mathbb{T}_{2\pi})}\leq C\epsilon.
\ee
}\end{theorem}

\br\label{subsonic-regularity}
{\it In fact, in the uniformly subsonic flow region $\Omega_{us}$, the regularity of the transonic flows can be improved to be $(U_1,U_2, B, A)\in (C^{3,\alpha}(\overline{\Omega_{us}}))^2\times(C^{4,\alpha}(\overline{\Omega_{us}}))^2$, where $\Omega_{us}=\{(r,\theta): r_c+\eta_0< r< r_1, \theta\in \mathbb{T}_{2\pi}\}$ for any $0<\eta_0<r_1-r_c$.}
\er

\br\label{nonzerovorticity}
{\it
Compared with the existence results of continuous subsonic-sonic or smooth transonic flows obtained in \cite{WX2013,WX2019,WX2020}, the flow constructed in Theorem \ref{2dmain} can have a small nonzero vorticity. As far as we know, this is the first result about the existence of nontrivial smooth transonic flows with nonzero vorticity.
}\er

\begin{remark}\label{mach}
{\it  In the case $|{\bf M}_b(r_0)|^2>\frac4{3-\gamma}$, since the total Mach number $|{\bf M}_b(r)|^2$ is strictly decreasing as $r$ increases, there exists a $\tilde{r}_0\in (r_0, r_1)$ such that $\frac4{3-\gamma}>|{\bf M}_b(\tilde{r}_0)|^2>1>M_{b1}^2(\tilde{r}_0)$. One may focus on the transonic flow region $\tilde{\Omega}:=\{\tilde{r}_0<r<r_1, \theta\in [0,2\pi]\}$. The extension of flow from $\tilde{\Omega}$ to $\Omega$ can be obtained by the well-developed theory for hyperbolic equations since the flow is purely supersonic in the subregion $\{(r,\theta): r_0<r<\tilde{r}_0\}$. The restriction \eqref{incoming} is acceptable in this sense. If one considers further the case with $\gamma\geq 3$, the condition \eqref{incoming} is automatically satisfied.
}\end{remark}

\begin{remark}\label{nonzeroradialvelocity}
{\it The requirement that the background flow should have a nonzero radial velocity can be removed if one considers only perturbations within the class of irrotational flows, i.e. $B_1(\theta)= A_1(\theta)\equiv 0$. In that case, there is no hyperbolicity in subsonic regions. Furthermore, the regularity assumptions on the boundary datum are weaken to be $(g_0,g_1)\in (H^3(\Omega))^2$. See Theorem \ref{main1} for more details.
}\end{remark}

\begin{remark}\label{higherregularity}
{\it There is a loss of derivatives in estimating the mixed type first order partial differential equations. To recover it, we introduce the stream function and observe that the regularity of the transonic flows in the uniformly subsonic region can be improved if the boundary datum at the entrance have better regularity. See Section \ref{proofTh12} for detailed explanations.
}\end{remark}

\begin{remark}\label{decelerating}
{\it There exists also a class of smooth irrotational transonic flows moving from the inner cylinder to the outer one (see \cite{WXY20a}). The speed of these flows decelerates smoothly from supersonic to subsonic. Such a transonic flow is also structurally stable under the same perturbations as in \eqref{2dboundary3}-\eqref{2dboundary4} within the class of irrotational flows. However, it seems difficult to improve the regularity of the flow in the uniformly supersonic region near the entrance. Therefore the loss of derivatives can not be recovered by the way developed in Section \ref{proofTh12} and it is not clear whether Theorem \ref{2dmain} holds or not for this kind of transonic flows.
}\end{remark}

Next we examine the structural stability of the background flow under suitable axisymmetric perturbations, therefore the problem simplifies to find axi-symmetric transonic flows in a concentric cylinder $\mathbb{D}=(r_0,r_1)\times \mathbb{R}$ satisfying suitable boundary conditions on the inner and outer cylinder respectively. Assume that the velocity and the density are of the form
\begin{eqnarray}\label{cyl1}
{\bf u}(x)= U_1(r,x_3){\bf e}_r + U_{2}(r,x_3) {\bf e}_{\theta} + U_3(r,x_3){\bf e}_3, \ \ \rho(x)=\rho(r,x_3),\ \ A(x)=A(r,x_3),
\end{eqnarray}
then the system \eqref{ProblemImm} reduces to
\begin{eqnarray}\label{cyl-sym}
\begin{cases}
\partial_r(\rho U_1)+\frac{1}{r} \rho U_1+ \partial_3(\rho U_3)=0,\\
(U_1\partial_r +U_3 \partial_3) U_1+\frac{1}{\rho} \partial_r p-\frac{U_2^2}r=0,\\
(U_1\partial_r +U_3 \partial_3) U_2+\frac{U_1U_2}r=0,\\
(U_1\partial_r +U_3 \partial_3) U_3+\frac{1}{\rho} \partial_{3}p=0,\\
(U_1\partial_r +U_3 \partial_3) A=0.
\end{cases}
\end{eqnarray}


The system \eqref{cyl-sym} is an elliptic-hyperbolic coupled system, in which $rU_2, A$ and $B$ satisfy the transport equations. So it is reasonable to prescribe the following boundary conditions on the outer cylinder $\{(r_{1},x_3): x_3\in\mathbb{R}\}$:
\begin{eqnarray}\label{cond1}
\begin{cases}
U_{2}(r_1,x_3)=U_{20}+\epsilon q_{2}(x_3),\  U_{3}(r_1,x_3)=\epsilon q_{3}(x_3)\\
A(r_1,x_3)=A_{1}+\epsilon \tilde{A}_1(x_3),\ B(r_1,x_3)=B_{1}+\epsilon \tilde{B}_1(x_3),
\end{cases}
\end{eqnarray}
with $q_{1}(x_3)$, $q_{2}(x_3)$, $\tilde{A}_1(x_3)$, $\tilde{B}_1(x_3)\in C^{2,\alpha}(\mathbb{R})$ and $\epsilon$ small enough. The boundary condition posed on the inner cylinder $\{(r_0,x_3): x_3\in\mathbb{R}\}$ is
\begin{eqnarray}\label{cond2}
U_{1}(r_0,x_3)=U_{b1}(r_0)+\epsilon q_{1}(x_3),
\end{eqnarray}
with $q_{3}(x_3)\in C^{2,\alpha}(\mathbb{R})$. For simplicity, we also assume that $q_k,k=1,2,3$ and $\tilde{B}_1, \tilde{A}_1$ have compact supports. It is expected that the flow will tend to the background state as $x_3\to \pm\infty$:
\be\label{cond3}
\displaystyle\lim_{x_3\to \pm\infty} (U_1,U_2,U_3,P,A)(r,x_3)= (U_{b1}(r),U_{b2}(r),0,P_0, A_0).
\ee

Then the following theorem on the existence and uniqueness of smooth axi-symmetric transonic spiral flows with small nonzero vorticity holds.

\begin{theorem}\label{cylind-theorem1}
{\it Given any background flow with nonzero radial velocity $U_{b1}\neq 0$, for any smooth $C^{2,\alpha}(\mathbb{R})$ functions $\tilde{B}_1, \tilde{A}_1$ and $q_k, k=1,2,3$ with compact supports, there exists a small constant $\epsilon_0$ depending only the background flow and boundary datum, such that if $0<\epsilon\leq \epsilon_0$, there exists a unique smooth transonic flow with nonzero vorticity
$${\bf u}=U_1(r,x_3){\bf e}_r + U_{2}(r,x_3) {\bf e}_{\theta} + U_3(r,x_3){\bf e}_3, \ \ A(x)=A(r,x_3), \ \ B(x)=B(r,x_3)$$
to \eqref{cyl-sym} with \eqref{cond1},\eqref{cond2} and \eqref{cond3}, and the following estimate holds
\begin{eqnarray}\label{cyl3}
\sum_{j=1}^2\|U_j-U_{bj}\|_{C^{2,\alpha}(\mathbb{D})}+ \|U_3\|_{C^{2,\alpha}(\mathbb{D})}+\|B- B_{0}\|_{C^{2,\alpha}(\mathbb{D})}+ \|A- A_{0}\|_{C^{2,\alpha}(\mathbb{D})}\leq C\epsilon,
\end{eqnarray}
for some constant $C$ depending only on the background solution and the boundary datum.

Moreover, all the sonic points form an axisymmetric surface with a parametric representation $r=\chi(x_3)\in C^1(\mathbb{R})$ extending from $-\infty$ to $\infty$. The sonic surface is closed to the background sonic cylinder in the sense that
\be\label{cy14}
\|\chi(x_3)-r_c\|_{C^1(\mathbb{R})}\leq C\epsilon
\ee
and
\be\label{cy15}
\displaystyle\lim_{x_3\to\pm\infty} \chi(x_3)= r_c.
\ee
}\end{theorem}

\br\label{remark4}
{\it Similar results can be obtained if one prescribes the vertical velocity $U_3$ and the pressure $P$ at the inner and outer cylinder respectively, instead of the radial velocity $U_1$ and the vertical velocity $U_3$ in \eqref{cond1}-\eqref{cond2}.
}\er

\br\label{cs-remark1}
{\it The requirement that the background flow has a nonzero radial velocity can be removed for perturbations within the classes of axi-symmetric transonic irrotational flows, i.e. $q_2= \tilde{A}_1=\tilde{B}_1\equiv 0$.
}\er

\br\label{po}
{\it There may exist another class of smooth irrotational transonic flow ${\bf u}(x)= \p_r \phi(r,x_3) {\bf e}_r + \frac{\kappa_2}{r}{\bf e}_{\theta}+ \p_{x_3}\phi(r,x_3) {\bf e}_3$ with a nonzero constant $\kappa_2$, which is not adjacent to the background transonic flow. The potential function $\phi$ will satisfy a second order mixed type differential equation with coefficients depending not only on $|\nabla \phi|^2$ but also on the space variable $r$. Some new difficulties arise from applying the Bernstein's method to get a fine gradient estimate. This will be reported in a forthcoming paper.
}\er

The theory of transonic fluid flows is closely related to the studies of the well-posedness theory for the mixed type partial differential equations. There are several classical mixed type PDEs which are closely related to the transonic fluid flows, such as Tricomi's equation, Keldysh's equation and the Von Karman equation, one may refer to \cite{Morawetz1958,Morawetz1970,OD1955,Protter1954,Smirnov1978} and the references therein for more details. Morawetz \cite{Morawetz1956,Morawetz1964} proved the nonexistence of a smooth solution to the perturbation for flow with a local supersonic region over a solid airfoil. Friedrichs \cite{Friedrichs1958} initiated a general and powerful theory of positive symmetric systems of first order and there are many important further development and applications to boundary value problems of equations of mixed type \cite{Gu1981,LP1969,LMP07,Morawetz2004}. Kuzmin \cite{Kuzmin2002} had investigated the nonlinear perturbation problem of an accelerating smooth transonic irrotational basic flow with some artificial boundary conditions in the potential and stream function plane. However, the existence of such a basic flow to the Chaplygin equation was not shown and the physical meaning of the boundary conditions was also not clear.

The existence of subsonic-sonic weak solutions to the 2-D steady potential equation were proved in \cite{cdsw07,xx07,xx10} by utilizing the compensated compactness and later on the authors \cite{chw16,hww11} examined the subsonic-sonic limit for multidimensional potential flows and steady Euler flows. However, the solutions obtained by the subsonic-sonic limit only satisfy the equations in the sense of distribution and there is no information about the regularity and degeneracy properties near sonic points and their distribution in flow region. Recently, in a series of papers \cite{WX2013,WX2016,WX2019,WX2020}, Wang and Xin have established the existence and uniqueness of Lipschitz continuous subsonic-sonic flows and smooth transonic flows of Meyer type in De Laval nozzles with a detailed description of sonic curve. In particular, under the assumption that the nozzle is suitably flat at its throat, they showed the existence and uniqueness of smooth transonic irrotational flows of Meyer type. The sonic points can locate only at the throat of the nozzle and the points on the nozzle wall with positive curvature. The sonic points at the throat are exceptional and strongly degenerate in the sense that all the characteristics from sonic points coincide with the sonic line and can not approach the supersonic region.

We make some comments on the key ingredients in our mathematical analysis for Theorem \ref{2dmain} and \ref{cylind-theorem1}. The authors in \cite{WX2013,WX2016,WX2019,WX2020} employed the Chaplygin equations in the plane of the velocity potential and the stream function and used the comparison principle as a main tool to analyze the subsonic-sonic flows. However, due to the nonzero angular velocity, the sonic points in our case are expected to be nonexceptional and the transformed sonic curve in the potential-stream functions plane is not a straight line in general, which is different from the cases studied by Wang and Xin, it seems quite difficult to adapt their methods to our case.  We need to find a different approach to deal with the flow with nontrivial vorticity. The steady Euler system is elliptic-hyperbolic mixed in subsonic regions and degenerates at sonic points. To circumvent this obstacle, we utilize the deformation-curl decomposition for the steady Euler system established by the first two authors in \cite{WengXin19,weng2019} to effectively decouple the hyperbolic and elliptic modes. This decomposition is based on a simple observation that one can rewrite the density equation as a Frobenius inner product of a symmetric matrix and the deformation matrix by using the Bernoulli's law. The vorticity is resolved by an algebraic equation of the Bernoulli's function and the entropy.

To explain the key ideas clearly, we first investigate the well-posedness theory to the linearized mixed type second order equation within the class of irrotational flows. By exploring some key properties of the background flows, we are able to find a class of multipliers and identify a class of admissible boundary conditions for the linearized problem, and this helps to yield the basic energy estimate and the high order derivatives estimates. Galerkin's method with Fourier series will be used for the construction of the approximated solutions and a simple contraction mapping argument will yield the solution to the nonlinear problem.

To further treat the rotational flows, note that the basic energy estimate for the linearized mixed type potential equation only helps to gain one order derivative regularity (see the $H^1$ estimate in Lemma \ref{basic-energy}), so that the iteration designed by the first two authors in \cite{WengXin19} for purely subsonic flows does not work in this case. We will choose some appropriate function spaces to design an elaborate two-layer iteration scheme to find the fixed point to the nonlinear problem. By requiring one order higher regularity of the boundary datum for the Bernoull's function and the entropy than those of the flow angles, we gain one more order derivatives estimates for the Bernoulli's function and the entropy than the velocity with the help of the stream function and the higher regularity of the flows in subsonic region. This is crucial for us to close the energy estimates.

The analysis of axi-symmetric transonic flows turns out to be simpler than those of the cylindrical transonic spiral flows. Again using the deformation-curl decomposition for the steady Euler equations in \cite{WengXin19,weng2019}, it will be shown that $U_1$ and $U_3$ satisfy a first order elliptic system when linearized around the background transonic flows. The quantities $r U_2, B$ and $S$ are conserved along the particle trajectory. The maximum principle and some suitable barrier functions are employed to obtain some uniform estimates to the second order elliptic equation. The far field behavior will be examined by a blow-up argument.

This paper will be arranged as follows. In Section \ref{key}, we state some key properties of the background flows which will play an important role in searching for an appropriate multiplier to the linearized mixed type potential equation. In Section \ref{annulus}, we first establish the basic and higher order energy estimates to the linearized mixed potential equations and construct approximated solutions by a Galerkin method. Then we employ the deformation-curl decomposition for the steady Euler system and design a two-layer iteration to demonstrate the existence of smooth transonic rotational flows. In Section \ref{cylinder}, we consider the structural stability of the background flows within the class of axi-symmetric flows.

\section{Some properties of the background flow}\label{key}\noindent

In the following, we derive some special properties of the background flow, which plays a key role in establishing the basic energy estimate for the linearized mixed type potential equation. It follows from \eqref{radial-euler} that
\begin{eqnarray}\nonumber
\begin{cases}
\rho_b'(r)=\frac{(M_{b1}^2+ M_{b2}^2)}{r(1-M_{b1}^2)}\rho_b ,\\
U_{b1}'(r)=-\frac{1+ M_{b2}^2}{r(1-M_{b1}^2)}U_{b1},\ \ U_{b2}'(r)=-\frac1r U_{b2},\\
(M_{b1}^2)'(r)=\frac{M_{b1}^2}{r(M_{b1}^2-1)}\left(2+(\gamma-1)M_{b1}^2+(\gamma+1)M_{b2}^2\right),\\
(M_{b2}^2)'(r)=\frac{M_{b2}^2}{r(M_{b1}^2-1)}\left(2+(\gamma-3)M_{b1}^2+(\gamma-1)M_{b2}^2\right),\\
(|{\bf M}_{b}|^2)'(r)=\frac{|{\bf M}_{b}|^2}{r(M_{b1}^2-1)}\left(2+(\gamma-1)|{\bf M}_{b}|^2\right).
\end{cases}
\end{eqnarray}
Thus the total Mach number of the background flow monotonically increases as $r$ decreases since $M_{b1}^2(r)<1$ for any $r\in [r_0,r_1]$.


For later use, we define
\begin{eqnarray}\nonumber
&&e_1(r)=\frac{c^2(\rho_b)+U_{b2}^2}r +\frac{(\ga+1)(1+M_{b2}^2)}{r(1-M_{b1}^2)}U_{b1}^2-\frac{(\ga-1)U_{b1}^2}r>0, \ \forall r\in [r_0,r_1],\\\nonumber
&&e_2(r)=\frac{2(1-M_{b1}^2)+(\ga-1)|{\bf M}_b|^2}{1-M_{b1}^2}\frac{U_{b1}U_{b2}}{r^2} ,
\end{eqnarray}
and
\begin{eqnarray}\nonumber
&&f(r)=\int_{r_0}^r\frac{M_{b1}M_{b2}(\tau)}{1-M_{b1}^2(\tau)}\frac{d\tau}\tau, \ \ k_{b33}(r)=\frac{1}{1-M_{b1}^2(r)},\\\nonumber
&&k_{b22}(r)=\frac{1-|{\bf M}_b(r)|^2}{r^2(1-M_{b1}^2(r))^2},\\\nonumber
&&k_{b1}(r)=\frac{e_1}{c^2(\rho_b)-U_{b1}^2}=\frac{1+M_{b2}^2+2M_{b1}^2+(\gamma+1)M_{b1}^2|{\bf M}_b|^2/(1-M_{b1}^2)}{r(1-M_{b1}^2)},\\\nonumber
&&k_{b2}(r)=f''(r)+\frac{e_1(r) f'(r)+ e_2(r)}{c^2(\rho_b)-U_{b1}^2}.
\end{eqnarray}

\begin{proposition}\label{key-property}
{\it Let $(U_{b1}, U_{b2}, \rho_b, A_0)$ be the background transonic flow, then the following identities hold for any $r\in [r_0,r_1]$
\begin{eqnarray}\label{key1}
&& k_{b2}(r)\equiv 0, \\\label{key2}
&& 2k_{b1}k_{b22}+k_{b22}'(r)=\frac{M_{b1}^2+M_{b2}^2}{r^3(1-M_{b1}^2)^3}\bigg(4-(3-\gamma)(M_{b1}^2+M_{b2}^2)\bigg),\\\label{key3}
&& 2k_{b1}k_{b33}+k_{b33}'(r)=\frac{2+2M_{b2}^2+(\gamma-1)M_{b1}^2 |{\bf M}_b|^2}{r(1-M_{b1}^2)^3}>0.
\end{eqnarray}
}\end{proposition}

\begin{proof}
To simplify notations, we denote $A_{b11}(r)=c^2(\rho_b)-U_{b1}^2$ and $A_{b12}(r)=-\frac{1}{r}U_{b1}U_{b2}$. Then $f'(r)=-\frac{A_{b12}}{A_{b11}}$. Direct calculations show that $k_{b2}=\frac{I}{A_{b11}^2}$, where
\begin{eqnarray}\nonumber
I&=&-e_1A_{b12}+ A_{b11} e_2- A_{b11} A_{b12}'+ A_{b12} A_{b11}'\\\nonumber
&=&\frac{U_{b1}U_{b2}}{r}\left\{\frac{c^2(\rho_b)+U_{b2}^2}r +\frac{(\ga+1)(1+M_{b2}^2)}{r(1-M_{b1}^2)}U_{b1}^2-\frac{(\ga-1)U_{b1}^2}r\right\}\\\nonumber
&\quad&+ A_{b11}\left\{\frac{2(1-M_{b1}^2)+(\ga-1)|{\bf M}_b|^2}{1-M_{b1}^2}\frac{U_{b1}U_{b2}}{r^2}-\frac{U_{b1}U_{b2}}{r^2}+ \frac{U_{b1}' U_{b2}+  U_{b1} U_{b2}'}{r}\right\}\\\nonumber
&\quad&-A_{b12}\{(\gamma+1)  U_{b1} U_{b1}'+(\gamma-1) U_{b2} U_{b2}'\}\\\nonumber
&=& \frac{(2-\gamma)U_{b1} U_{b2}}{r^2}\{U_{b1}^2 + U_{b2}^2 - (c^2(\rho_b)+U_{b2}^2)+A_{b11}\}\\\nonumber
&\equiv& 0.
\end{eqnarray}

Moreover, one can calculate that
\begin{eqnarray}\nonumber
&&k_{b22}'(r)=\frac{-2(1-|{\bf M}_b|^2)}{r^3(1-M_{b1}^2)^2}-\frac{2(1-|{\bf M}_b|^2)}{r^2(1-M_{b1}^2)^3} \frac{M_{b1}^2[2(1-M_{b1}^2)+(\gamma+1)|{\bf M}_b|^2]}{r(1-M_{b1}^2)}\\\no
&&\quad\quad\quad\quad+ \frac{1}{r^2(1-M_{b1}^2)^2}\frac{|{\bf M}_b|^2(2+(\gamma-1)|{\bf M}_b|^2)}{r(1-M_{b1}^2)}\\\no
&&=\frac{(\gamma-1)|{\bf M}_b|^4+4 |{\bf M}_b|^2-2- 2M_{b1}^2(1-|{\bf M}_b|^2)}{r^3(1-M_{b1}^2)^3}-\frac{2(\gamma+1)M_{b1}^2|{\bf M}_b|^2(1-|{\bf M}_b|^2)}{r^3(1-M_{b1}^2)^4}.
\end{eqnarray}

Therefore
\begin{eqnarray}\nonumber
&&2k_{b1}k_{b22}+k_{b22}'=\frac{2(1-|{\bf M}_{b}|^2)}{r^3(1-M_{b1}^2)^3}\left(1+M_{b1}^2+|{\bf M}_{b}|^2+\frac{(\gamma+1)M_{b1}^2|{\bf M}_{b}|^2}{1-M_{b1}^2}\right)\\\nonumber
&&+\frac{(\gamma-1)|{\bf M}_b|^4+4 |{\bf M}_b|^2-2- 2M_{b1}^2(1-|{\bf M}_b|^2)}{r^3(1-M_{b1}^2)^3}-\frac{2(\gamma+1)M_{b1}^2|{\bf M}_b|^2(1-|{\bf M}_b|^2)}{r^3(1-M_{b1}^2)^4}\\\nonumber
&&=\frac{|{\bf M}_{b}|^2}{r^3(1-M_{b1}^2)^3}\bigg(4-(3-\gamma)|{\bf M}_{b}|^2\bigg).
\end{eqnarray}


On the other hand, we have
\begin{eqnarray}\nonumber
k_{b33}'(r)= -\frac{M_{b1}^2}{r(1-M_{b1}^2)^3}\left\{2+(\gamma-1)M_{b1}^2+(\gamma+1)M_{b2}^2\right\}
\end{eqnarray}
and
\begin{eqnarray}\nonumber
&&2k_{b1}k_{b33}+k_{b33}'= \frac{2}{r(1-M_{b1}^2)^3}\left\{(1+|{\bf M}_b|^2+M_{b1}^2)(1-M_{b1}^2)+(\gamma+1)M_{b1}^2|{\bf M}_b|^2\right\}\\\nonumber
&&\quad\quad -\frac{M_{b1}^2}{r(1-M_{b1}^2)^3}\left\{2+(\gamma-1)M_{b1}^2+(\gamma+1)M_{b2}^2\right\}\\\nonumber
&&=\frac{2+2M_{b2}^2+(\gamma-1)M_{b1}^2 |{\bf M}_b|^2}{r(1-M_{b1}^2)^3}.
\end{eqnarray}
Thus the proposition is proved.

\end{proof}

\section{Smooth cylindrical transonic flows with nonzero vorticity}\label{annulus}\noindent

Since the background flow changes smoothly from subsonic at the outer circular cylinder to supersonic at the inner one, the linearized potential equation is of mixed type in $\Omega$. We would concentrate on searching for an appropriate multiplier and identifying suitable boundary conditions for the linearized mixed type potential equation. To illustrate the main ideas, we start with the potential flows.

\subsection{Irrotational flows}\noindent

For a potential flow, the vorticity being free implies that, in terms of the polar coordinates,
\begin{eqnarray}\label{irrotational}
\frac{1}{r}\p_r(r U_2)-\frac{1}{r}\p_\th U_1=0,\ \ \text{in}\ \Omega
\end{eqnarray}
and
\begin{eqnarray}\label{density}
\rho=\left(\frac{\gamma-1}{A_0\gamma}\right)^{\frac{1}{\gamma-1}}\left(B_0-\frac{1}2(U_1^2+U_2^2)\right)^{\frac{1}{\gamma-1}}.
\end{eqnarray}
Therefore \eqref{2d} is reduced to the following boundary value problem in $\Omega$:
\begin{align}\label{2dcom-euler}
\begin{cases}
\partial_r(r\rho U_1)+\partial_\theta(\rho U_2)=0,\\
\frac{1}{r}\partial_r(r U_2)-\frac{1}{r}\partial_\theta U_1=0,\\
U_1(r_0,\theta)-l_0 U_2(r_0,\theta)=a_0+\epsilon g_0(\theta),\\
U_2(r_1,\theta)=a_1+\epsilon g_1(\theta).
\end{cases}
\end{align}
Recall that $l_0$ is a constant in a suitable range to be specified and
\begin{eqnarray}\nonumber
a_0= U_{b1}(r_0)- l_0 U_{b2}(r_0),\ \ \ a_1=U_{20}.
\end{eqnarray}

The existence and uniqueness of smooth irrotational transonic flows can written as follows.
\begin{theorem}\label{main1}
{\it
Let the background flow and $l_0$ be given as in Theorem \ref{2dmain} except the assumption $U_{b1}\neq 0$. Assume that $g_0, g_1\in H^3(\mathbb{T}_{2\pi})$. Then there exists a small constant $\epsilon_0$ depending on the background flow, $l_0$ and $g_0, g_1$, such that for any $0<\epsilon<\epsilon_0$, the problem \eqref{2dcom-euler} has a unique smooth transonic irrotational solution $(U_1,U_2)\in H^3(\Omega)\subset C^{1,\alpha}(\bar{\Omega})$ with the estimate
\be\label{estimate1}
\|U_1-U_{b1}\|_3+\|U_2-U_{b2}\|_3\leq C\epsilon.
\ee
Moreover, all the sonic points form a closed arc with a parametric representation $r=s(\th)\in C^{1}(\mathbb{T}_{2\pi})$ with any $\alpha\in (0,1)$. The sonic curve is closed to the background sonic circle in the sense that
\be\label{estimate2}
\|s(\theta)-r_c\|_{C^1(\mathbb{T}_{2\pi})}\leq C\epsilon.
\ee
}\end{theorem}

\br\label{remark2}
{\it Compared with Theorem \ref{2dmain}, the assumption that the background solution should have a nonzero radial velocity $U_{b1}\neq 0$ is removed for potential flows here, since there is no need to solve the transport equations for the Bernoulli's function and the entropy in this case.
}\er

\br\label{remark3}
{\it All the sonic points to the transonic irrotational solution obtained in Theorem \ref{main1} are nonexceptional and noncharacteristically degenerate.
}\er

\br\no
{\it The cylindrically symmetric smooth transonic flows where the fluid moves from the inner to the outer circle are also structurally stable under the same perturbations as in \eqref{2dcom-euler} within the class of irrotational flows. 
}\er

Since $\Omega$ is non simply connected and the background flow has a nonzero circulation, the potential function corresponding to the background flow is $\phi_b(r,\theta)=\int_{r_1}^r U_{b1}(s)ds + r_1 U_{20} \theta$, which is not periodic in $\theta$. To avoid the trouble, we denote the difference between the flow and the background flow by
\begin{eqnarray}\label{difference}
\hat{U}_1= U_1- U_{b1},\ \ \hat{U}_2= U_2- U_{b2},\ \ \hat{\rho}= \rho-\rho_b,
\end{eqnarray}
then ${\bf \hat{U}}$ and $\hat{\rho}$ satisfy
\begin{eqnarray}\label{differ1}\begin{cases}
\partial_r\left(r(\rho_b \hat{U}_1+ (U_{b1}+\hat{U}_1)\hat{\rho})\right)+ \partial_{\theta}\left(\rho_b \hat{U}_2+ (U_{b2}+\hat{U}_2)\hat{\rho}\right)=0,\\
\frac{1}{r}\partial_{r} (r \hat{U}_2)- \frac{1}{r}\partial_{\theta} \hat{U}_1=0,\\
\hat{U}_1(r_0,\theta)-l_0 \hat{U}_2(r_0,\theta)=\epsilon g_0(\theta),\\
\hat{U}_2(r_1,\theta)=\epsilon g_1(\theta).
\end{cases}
\end{eqnarray}

Define the potential function
\begin{align}\no
\phi(r,\theta)&=\int_{r_1}^r\hat{U}_1(\tau,\theta)d\tau+\int_0^\theta(r_1\hat{U}_2(r_1,\tau)+d_0)d\tau,\\\no
&=\int_{r_1}^r \hat{U}_1(\tau,\th)d\tau+\int_0^{\theta}(\epsilon r_1 g_1(\tau)+ d_0) d\tau,
\end{align}
where $d_0$ is introduced so that $\phi(r,\theta)=\phi(r,\theta+2\pi)$. Indeed, $d_0=-\epsilon r_1\bar{g}_1$ with $\bar{g}_1\equiv\frac1{2\pi}\int_0^{2\pi}g_1(\tau)d\tau$. Then $\phi$ is periodic in $\theta$ with period $2\pi$ and satisfy
\begin{align}\label{phi00}
\partial_r\phi=\hat{U}_1,\ \ \partial_\theta\phi=r \hat{U}_2+d_0.
\end{align}

Substituting \eqref{phi00} into \eqref{differ1} yields that $\phi$ satisfies a second order mixed type equation
\begin{align}\label{Pro}
\begin{cases}
L\phi\equiv A_{11}\partial_r^2\phi+A_{22}\partial_\th^2\phi+(A_{12}+A_{21})\partial_{r\th}^2\phi+e_1(r)\partial_r\phi+e_2(r)\partial_\th\phi=F(U_1,U_2),\\
\partial_r\phi(r_0,\th)-l_0\frac1{r_0}\partial_\th\phi(r_0,\th)=\epsilon g_0(\th)-l_0\frac1{r_0} d_0,\\
\partial_\th\phi(r_1,\th)=d_0+r_1\epsilon g_1(\th),\\
\phi(r_1,0)=0,
\end{cases}
\end{align}
where
\begin{eqnarray}\nonumber
A_{11}(U_1,U_2)&=& c^2(\rho)-U_1^2,\  \  A_{22}(U_1,U_2)=\frac{c^2(\rho)- U_2^2}{r^2},\\\nonumber
A_{12}(U_1,U_2)&=& A_{21}(U_1,U_2)=-\frac{U_1 U_2}{r},\\\nonumber
F(U_1,U_2)&=& e_2(r)d_0+(\frac{\gamma+1}2 U_{b1}'+\frac{\gamma-1}{2r}U_{b1})\hat{U}_1^2+(\frac{\gamma-1}2 U_{b1}'+\frac{\gamma-3}{2r}U_{b1})\hat{U}_2^2\\\nonumber
&\quad&-\frac{c^2(\rho)+U_2^2-c^2(\rho_b)-U_{b2}^2}r \hat{U}_1.
\end{eqnarray}

\subsubsection{Linearized problem}\noindent

Denote the function space
\begin{eqnarray}\nonumber
\mathcal{X}=\{\phi\in H^4(\Omega), \|\phi\|_4\leq \delta_0\},
\end{eqnarray}
where $\delta_0>0$ will be specified later. For any function $\bar{\phi}\in \mathcal{X}$, define
\begin{eqnarray}\nonumber
\bar{U}_1= U_{b1}+\partial_r \bar{\phi},\ \ \bar{U}_2= U_{b2}+ \frac{1}{r}\partial_{\theta} \bar{\phi}- \frac{d_0}{r}.
\end{eqnarray}

We will construct an operator $\mathcal{T}$: $\bar{\phi}\in \mathcal{X}\mapsto\phi\in \mathcal{X}$, where $\phi$ will be obtained by solving the following linear mixed-type second-order partial differential equation
\begin{align}\label{LinearPotential2}
\begin{cases}
\bar{L}\phi\equiv A_{11}(\bar{U}_1,\bar{U}_2)\p_r^2\phi+A_{22}(\bar{U}_1,\bar{U}_2)\p_\th^2\phi+2A_{12}(\bar{U}_1,\bar{U}_2)\p_{r\theta}^2\phi\\
\ \ \ \ \ +e_1(r)\p_r\phi
+e_2(r)\p_\th\phi=F(\bar{U}_1,\bar{U}_2),\\
\partial_r\phi(r_0,\th)-\frac{l_0}{r_0}\partial_\th\phi(r_0,\th)=\epsilon g_0(\th)-\frac{l_0}{r_0}d_0=O(\eps),\\
\partial_\th\phi(r_1,\th)=d_0+r_1\epsilon g_1(\th)=O(\eps),\\
\phi(r_1,0)=0,
\end{cases}
\end{align}
with the coefficients satisfying the following estimates
\begin{align}
\begin{cases}
\|A_{ij}(\bar{U}_1,\bar{U}_2)-A_{ij}(U_{b1},U_{b2})\|_3\leq C_0\delta_0,\ i,j=1,2,\\
\|F(\bar{U}_1,\bar{U}_2)\|_3\leq C_0(\epsilon+\delta_0^2).
\end{cases}
\end{align}
Here and in the following the constant $C_0$ depends only the background flow, the boundary datum and $l_0$ and may change from line to line.

Define a new coordinate $(y_1, y_2)$ as
\begin{eqnarray}\nonumber
y_1=r,\ \ y_2= f(r)+\theta,\
\end{eqnarray}
where
\begin{eqnarray}\nonumber
f(r)&=&-\int_{r_0}^r\frac{A_{b12}(\tau)}{A_{b11}(\tau)}d\tau
=\int_{r_0}^r\frac{M_{b1}M_{b2}(\tau)}{1-M_{b1}^2(\tau)}\frac{d\tau}\tau.
\end{eqnarray}
Set the function $\hat{\phi}(y_1,y_2)= \phi(y_1,y_2-f(y_1))$. Then \eqref{LinearPotential2} can be written as
\begin{align}\label{LinearPotential3}
\begin{cases}
\hat{L}\hat{\phi}\equiv\displaystyle \sum_{i,j=1}^2 k_{ij}\partial_{y_iy_j}^2 \hat{\phi} + \sum_{i=1}^2 k_i\partial_{y_i}\hat{\phi}=\frac{F(\bar{U}_1,\bar{U}_2)}{A_{11}(\bar{U}_1,\bar{U}_2)}=:\hat{F}(\bar{U}_1,\bar{U}_2),\\
r_0\partial_{y_1}\hat{\phi}(r_0,y_2)+(r_0f'(r_0)-l_0)\partial_{y_2}\hat{\phi}(r_0,y_2)=g_2(y_2),\\
\partial_{y_2}\hat{\phi}(r_1,y_2)=g_3(y_2),\\
\hat{\phi}(r_1,f(r_1))=0,
\end{cases}
\end{align}
where $(y_1,y_2)\in (r_0,r_1)\times \mathbb{T}_{2\pi} $ and
\begin{eqnarray}\no
&&k_{11}\equiv 1,\ \ k_{12}(\bar{U}_1,\bar{U}_2)=k_{21}(\bar{U}_1,\bar{U}_2)= \frac{A_{12}(\bar{U}_1,\bar{U}_2)+A_{11}(\bar{U}_1,\bar{U}_2) f'(y_1)}{A_{11}(\bar{U}_1,\bar{U}_2)},\\\nonumber
&&k_{22}(\bar{U}_1,\bar{U}_2)=\frac{A_{22}(\bar{U}_1,\bar{U}_2)+(A_{12}+A_{21})(\bar{U}_1,\bar{U}_2)f'(y_1)}{A_{11}(\bar{U}_1,\bar{U}_2)}+(f'(y_1))^2,\\\nonumber
&&k_{1}(\bar{U}_1,\bar{U}_2)=\frac{e_1(y_1)}{A_{11}(\bar{U}_1,\bar{U}_2)},\ \ k_{2}(\bar{U}_1,\bar{U}_2)=f''(y_1)+\frac{e_1(y_1) f'(y_1)+ e_2(y_1)}{A_{b11}(\bar{U}_1,\bar{U}_2)},\\\no
&& g_2(y_2)=r_0\epsilon g_0(y_2)-l_0d_0,\ \ g_3(y_2)=d_0+r_1\epsilon g_1(y_2-f(r_1)).
\end{eqnarray}

 Then it follows from the definitions in Section \ref{key} and the facts that $A_{b12}(y_1)+A_{b11}(y_1)f'(y_1)\equiv0$ and $k_{b2}(y_1)\equiv 0$ in Proposition \ref{key-property} that the following important estimates hold
\begin{align}\label{coe}
\begin{cases}
\|k_{12}(\bar{U}_1,\bar{U}_2)\|_3+\|k_{22}(\bar{U}_1,\bar{U}_2)-k_{b22}(y_1)\|_3\leq C_0\delta_0,\\
\|k_1(\bar{U}_1,\bar{U}_2)-k_{b1}(y_1)\|_3+\|k_2(\bar{U}_1,\bar{U}_2)\|_3\leq C_0\delta_0.
\end{cases}
\end{align}

To simplify the notation, we still use $\phi$ instead of $\hat{\phi}$ in the following.

\subsubsection{Energy estimates for the linearized problem}\label{energy} \noindent

In this subsection, we will derive the energy estimate to \eqref{LinearPotential3} under the assumptions that $k_{ij}, k_i (i,j=1,2)\in C^{\infty}(\bar{\Omega})$ and $g_2, g_3\in C^{\infty}(\mathbb{T}_{2\pi})$ and \eqref{coe} holds.

\begin{lemma}\label{basic-energy}
{\it Suppose that \eqref{incoming} and \eqref{bo} hold. There exists two constants $\sigma_*>0, \delta_*>0$ depending only on the background flow, and $l_0$, such that if $0<\delta_0\leq \delta_*$ in \eqref{coe}, the solution to \eqref{LinearPotential3} satisfies the following basic energy estimate
\begin{eqnarray}\label{basic-energy1}
\|\phi\|_1\leq \frac{C_*}{\sigma_*} (\|\hat{F}\|_{L^2(\Omega)}+ \sum_{j=0,1} \|g_j\|_{L^2(\mathbb{T}_{2\pi})}),
\end{eqnarray}
where the constant $C_*$ depends only on the $H^3(\Omega)$ norms of the coefficients $k_{ij}, k_i$ for $i,j=1,2$.
}\end{lemma}

\begin{proof}
We employ the basic idea of positive operator theory developed by Friedrichs \cite{Friedrichs1958} and some key properties of the background flow to find a multiplier and identify a class of admissible boundary conditions at the inner and outer circle, which yields the basic energy estimate. Let $l_1(y_1)$ and $l_2(y_1)$ be smooth functions of $y_1$ in $[r_0,r_1]$ to be determined. Integration by parts leads to
\begin{align}\nonumber
&\iint_{\Omega}\hat{F}(l_1(y_1)\partial_{y_1}\phi+l_2(y_1)\partial_{y_2}\phi) dy_1dy_2\\\nonumber
=&\frac12\int_0^{2\pi}\left(\sqrt{l_1}\partial_{y_1}\phi+\frac{l_2}{\sqrt{l_1}}\partial_{y_2}\phi\right)^2+\left(-k_{22}l_1+2k_{12} l_2-\frac{l_2^2}{l_1}\right)(\partial_{y_2}\phi)^2 dy_2 \bigg|^{r_1}_{y_1=r_0}\\\no
&+\iint_{\Omega}\left(l_1k_1-\frac12 l'_1- l_1\partial_{y_2} k_{12}\right)(\partial_{y_1}\phi)^2+(k_1l_2-l_2'+l_1k_2-l_1\partial_{y_2}k_{22})\partial_{y_1}\phi\partial_{y_2}\phi\\\label{basic}
&\quad\quad\quad +\left(\frac12\partial_{y_1}(l_1k_{22})-\frac12 l_2\partial_{y_2}k_{22}-\partial_{y_1}(l_2k_{12})+l_2k_2\right)(\partial_{y_2}\phi)^2 dy_1dy_2.
\end{align}
To get an energy estimate, we will show that there exist smooth functions $l_1(y_1)$ and $l_2(y_1)$ such that if $\delta_0$ in \eqref{coe} is small enough, the following inequalities hold
\begin{eqnarray}\label{positive}\begin{cases}
l_1k_1-\frac12 l'_1- l_1\partial_{y_2} k_{12}\geq \sigma_*,\ \ \forall (x,y)\in \Omega,\\
\frac12\partial_{y_1}(l_1k_{22})-\frac12 l_2\partial_{y_2}k_{22}-\partial_{y_1}(l_2k_{12})+l_2k_2\geq \sigma_*,\ \ \forall (x,y)\in \Omega,\\
\|k_1l_2-l_2'+l_1k_2-l_1\partial_{y_2}k_{22}\|_{L^{\infty}(\Omega)}\leq C_*\delta_0,\\
\left(k_{22}l_1+\frac{l_2^2}{l_1}-2k_{12} l_2\right)(r_0,y_2)>0,\ \ \forall y_2\in\mathbb{T}_{2\pi},\\
\end{cases}\end{eqnarray}
where $\sigma_*>0$ is a constant depending only on the background flow.

Choose $l_1(y_1)$ such that $l_1(y_1)k_{b1}(y_1)-\frac12 l_1'(y_1)=\sigma_1>0$, where $\sigma_1\in (0, 1)$ is a small positive constant to be determined later. Therefore
\begin{align}\label{l2}
l_1(y_1)=e^{\int_{r_0}^{y_1} 2 k_{b1}(\tau_1)d\tau_1}\left(l_1(r_0)-2\sigma_1\int_{r_0}^{y_1} e^{-\int_{r_0}^{\tau_1}2k_{b1}(\tau_2)d\tau_2}d\tau_1\right),
\end{align}
where $l_1(r_0)=1+2\int_{r_0}^{r_1} e^{-\int_{r_0}^{\tau_1}2k_{b1}(\tau_2)d\tau_2}d\tau_1$. Then $l_1(y_1)>0$ for any $y_1\in [r_0,r_1]$ if $\sigma_1\leq \frac12$.

Recall the identity \eqref{key2} from Proposition \ref{key-property}, we have
\begin{eqnarray}\nonumber
(l_1k_{b22})'(y_1)&=&l_1'k_{b22}+l_1 k_{b22}'=l_1(2k_{b1}k_{b22}+k_{b22}')-2\sigma_1k_{b22}\\\nonumber
&=& \frac{|{\bf M_{b}}|^2}{y_1^3(1-M_{b1}^2)^3}\bigg(4-(3-\gamma)|{\bf M_{b}}|^2\bigg) l_1(y_1)- 2\sigma_1 \frac{1-|{\bf M}_b|^2}{y_1^2(1-M_{b1}^2)^2}.
\end{eqnarray}

Since $\frac{d}{dy_1} |{\bf M}_b|^2<0$ for any $r\geq r_0$, then if $|{\bf M}_{b}(r_0)|^2<\frac4{3-\gamma}$, it holds that $|{\bf M}_{b}(y_1)|^2<\frac4{3-\gamma}$ for any $y_1\in [r_0, r_1]$. Set $\sigma_0= \frac4{3-\gamma}-|{\bf M}_{b}(r_0)|^2$. There exist positive constants $\sigma_2$ and $\sigma_3$ depending only on the background flow such that if $0<\sigma_1\leq \sigma_2$, then
\begin{eqnarray}\nonumber
(l_1k_{b22})'(y_1)&\geq& \frac{1}{y_1^2(1-M_{b1}^2)^2}\left\{(3-\gamma)\sigma_0\frac{|{\bf M}_b|^2}{y_1(1-M_{b1}^2)}l_1(y_1)
-2\sigma_1(1-|{\bf M}_{b}|^2)\right\}\\\no
&\geq&\frac{1}{y_1^2(1-M_{b1}^2)^2}\bigg\{(3-\gamma)\sigma_0\frac{|{\bf M}_b|^2}{y_1(1-M_{b1}^2)}e^{\int_{r_0}^{y_1} 2 k_{b1}(\tau_1)d\tau_1}\bigg(l_1(r_0)\\\nonumber
&\quad&\quad-2\sigma_1\int_{r_0}^{y_1} e^{-\int_{r_0}^{\tau_1}2k_{b1}(\tau_2)d\tau_2}d\tau_1\bigg)-2\sigma_1(1-|{\bf M}_{b}|^2)\bigg\}\\\nonumber
&\geq& \sigma_3, \quad \forall y_1\in [r_0,r_1].
\end{eqnarray}

Set $\sigma_1=\sigma_2$ in \eqref{l2}. It follows from \eqref{coe} that for $\delta_0$ small enough, one has
\begin{eqnarray}\nonumber
&&l_1k_1-\frac12 l'_1- l_1\partial_{y_2} k_{12}=l_1 k_{b1}- \frac{1}{2} l_1'+ l_1(k_1-k_{b1})- l_1\partial_{y_2} k_{12}\\\label{coe3}
&&\geq \sigma_2- \|l_1(k_1-k_{b1})\|_{L^{\infty}}- \|l_1\partial_{y_2} k_{12}\|_{L^{\infty}}\geq \frac{1}{2}\sigma_2>0, \ \ \forall (x,y)\in \Omega,
\end{eqnarray}
due to the Sobolev embedding $H^3(\Omega)\subset C^{1,\alpha}(\overline{\Omega})$ with $\alpha\in(0,1)$.

Set
\begin{align}\no
l_2(y_1)=\left(f'(r_0)-\frac{l_0}{r_0}\right)l_1(r_0) e^{\int_{r_0}^{y_1} k_{b1}(\tau)d\tau},
\end{align}
where $l_0$ is a constant to be chosen such that
\begin{eqnarray}\label{coe1}
k_{b22}(r_0)+\left(f'(r_0)-\frac{l_0}{r_0}\right)^2>0,
\end{eqnarray}
which is equivalent to \eqref{bo} in Theorem \ref{2dmain}.

Therefore $l_2(y_1)$ satisfies $l_2(y_1)k_{b1}(y_1)-l_2'(y_1)=0$ and one should note that
\begin{eqnarray}\nonumber
\left(\sqrt{l_1}\partial_{y_1}\phi+\frac{l_2}{\sqrt{l_1}}\partial_{y_2}\phi\right)(r_0,y_2)
&=&\frac{\sqrt{l_1(r_0)}}{r_0}\left(r_0\p_{y_1}\phi+(r_0f'(r_0)-l_0)\p_{y_2}\phi\right)(r_0,y_2)\\\label{boundary2}
&=&\frac{\sqrt{l_1(r_0)}}{r_0} g_2(y_2).
\end{eqnarray}
Since $k_{b22}(r_1)>0$, if $\delta_0$ is small enough, then
\be\no
&\quad&\frac{l_2^2}{l_1}(r_1,y_2)+k_{22}(r_1,y_2)l_1(r_1)-2k_{12}(r_1,y_2)l_2(r_1)\\\no
&=&\frac{1}{l_1(r_1)}(l_2^2+k_{b22}l_1^2)(r_1)+(k_{22}-k_{b22})(r_1,y_2)l_1(r_1)- 2k_{12}(r_1,y_2)|l_2(r_1)|\\\no
&\geq&\frac{1}{l_1(r_1)}(l_2^2+k_{b22}l_1^2)(r_1)-\|k_{22}-k_{b22}\|_{L^{\infty}}l_1(r_1)- 2\|k_{12}\|_{L^{\infty}}|l_2(r_1)|\\\label{coe5}
&\geq& \frac{1}{2l_1(r_1)}(l_2^2+k_{b22}l_1^2)(r_1)>0
\ee
and
\be\no
&\quad&\frac{l_2^2}{l_1}(r_0,y_2)+k_{22}(r_0,y_2)l_1(r_0)-2k_{12}(r_0,y_2)l_2(r_0)\\\no
&=&l_1(r_0)\left(k_{b22}(r_0)+(f'(r_0)-\frac{l_0}{r_0})^2\right)+(k_{22}-k_{b22})(r_0,y_2)l_1(r_0)- 2k_{12}(r_0,y_2)l_2(r_0)\\\no
&\geq&l_1(r_0)\left(k_{b22}(r_0)+(f'(r_0)-\frac{l_0}{r_0})^2\right)-\|k_{22}-k_{b22}\|_{L^{\infty}}l_1(r_0)- 2\|k_{12}\|_{L^{\infty}}|l_2(r_0)|\\\label{coe6}
&\geq& \frac1{2} l_1(r_0)\left(k_{b22}(r_0)+(f'(r_0)-\frac{l_0}{r_0})^2\right)>0.
\ee

With $l_2(y_1)$ fixed, if $\delta_0$ is chosen small enough, then
\begin{eqnarray}\nonumber
&&\frac12\partial_{y_1}(l_1k_{22})-\frac12 l_2\partial_{y_2}k_{22}-\partial_{y_1}(l_2k_{12})+l_2k_2 \\\nonumber
&&=\frac{1}{2}(l_1 k_{b22})'(y_1)+\frac{1}{2}\partial_{y_1}(l_1(k_{22}-k_{b22}))-\frac12 l_2\partial_{y_2}k_{22}-\partial_{y_1}(l_2k_{12})+l_2k_2\\\nonumber
&&\geq \sigma_3-\|\partial_{y_1}(l_1(k_{22}-k_{b22}))\|_{L^{\infty}}-\|l_2\partial_{y_2}k_{22}\|_{L^{\infty}}-\|\partial_{y_1}(l_2k_{12})\|_{L^{\infty}}
-\|l_2k_2\|_{L^{\infty}}\\\no
&&\geq \frac{1}{2}\sigma_3,\ \ \ \forall (y_1,y_2)\in \Omega,\\\no
&&\|k_1l_2-l_2'+l_1k_2-l_1\partial_{y_2}k_{22}\|_{L^{\infty}}=\|l_2(k_1-k_{b1})-l_1 \partial_{y_2} k_{22}+ l_1 k_2\|_{L^{\infty}}\leq C_* \delta_0.
\end{eqnarray}
Hence the inequalities in \eqref{positive} are proved for $\sigma_*=\frac{1}{2}\min\{\sigma_2,\sigma_3\}$.

With the help of \eqref{positive} and \eqref{coe5}-\eqref{coe6}, one can conclude from \eqref{basic} that
\begin{eqnarray}\nonumber
&\quad&\iint_{\Omega}(|\partial_{y_1}\phi|^2+|\partial_{y_2}\phi|^2)dy_1dy_2\\\no
&\quad&\quad+\int_0^{2\pi}\left(\sqrt{l_1}\partial_{y_1}\phi+\frac{l_2}{\sqrt{l_1}}\partial_{y_2}\phi\right)^2(r_1,y_2)+ (\p_{y_2}\phi)^2(r_0,y_2)dy_2\\\no
&\leq& \frac{C_*}{\sigma_*}\bigg(\iint_{\Omega}|\hat{F}(y_1,y_2)|^2 dy_1dy_2
+\sum_{j=2}^3\int_0^{2\pi}|g_j(y_2)|^2dy_2\bigg).
\end{eqnarray}

Since $\phi(y_1,y_2)=\int_{r_1}^{y_1}\partial_{y_1}\phi(\tau, y_2)d\tau+\int_{f(r_1)}^{y_2}g_3(\tau)d\tau$, thus
\begin{align}\no
\|\phi\|_{1}^2\leq\frac{C_*}{\sigma_*}\bigg(\iint_{\Omega}\hat{F}^2(y_1,y_2)dy_1dy_2
+\sum_{j=2}^3\int_0^{2\pi}|g_j(y_2)|^2dy_2\bigg).
\end{align}
We have finished the proof of Lemma \ref{basic-energy}.
\end{proof}


\begin{lemma}\label{highorder}
{\it Under the assumptions of Lemma \ref{basic-energy}, the following high order derivatives estimate holds:
\begin{eqnarray}\label{highorder1}
\|\phi\|_4 \leq \frac{C_*}{\sigma_*} (\|\hat{F}\|_3+ \|g_0\|_{3} +\|g_1\|_3).
\end{eqnarray}
}\end{lemma}

\begin{proof}

Let $v=\partial_{y_2}\phi$. Then $v$ satisfies
\begin{eqnarray}\label{second1}
\begin{cases}
\displaystyle \sum_{i,j=1}^2 k_{ij} \partial_{y_i y_j}^2 v + (k_1+ 2 \partial_{y_2} k_{12})\partial_{y_1} v + (k_2+ \partial_{y_2} k_{22}) \partial_{y_2} v  \\
\quad\quad\quad\quad=\partial_{y_2}\hat{F}-\partial_{y_2}k_1\partial_{y_1}\phi-\partial_{y_2}k_2\partial_{y_2}\phi,\\
r_0\partial_{y_1} v(r_0,y_2)+(r_0f'(r_0)-l_0)\partial_{y_2} v(r_0,y_2)=g'_2(y_2),\\
\partial_{y_2} v(r_1, y_2)=g_3'(y_2),
\end{cases}
\end{eqnarray}
Multiplying the equation in \eqref{second1} by $l_1(y_1)\partial_{y_1} v+ l_2(y_1)\partial_{y_2} v$ and integrating over $\Omega$, after integration by parts, one can get
\begin{eqnarray}\nonumber
&&\iint_{\Omega}\left(l_1k_1-\frac12 l'_1+ l_1\partial_{y_2} k_{12}\right)|\partial_{y_1}v|^2+(k_1l_2-l_2'+k_2l_1+2l_2\partial_{y_2}k_{12})\partial_{y_1} v\partial_{y_2} v\\\label{second-energy}
&&+\left(\frac12\partial_{y_1}(l_1k_{22})+\frac12 l_2\partial_{y_2}k_{22}-\partial_{y_1}(l_2k_{12})+l_2k_2\right)|\partial_{y_2}v|^2 dy_1dy_2\\\no
&&+\frac12\int_0^{2\pi}\left(\sqrt{l_1}\partial_{y_1}v+\frac{l_2}{\sqrt{l_1}}\partial_{y_2}v\right)^2+\left(-k_{22}l_1+2k_{12} l_2-\frac{l_2^2}{l_1}\right)|\partial_{y_2}v|^2 dy_2 \bigg|^{r_1}_{y_1=r_0}\\\nonumber
&&=\iint_{\Omega}(\partial_{y_2}\hat{F}-\partial_{y_2} k_1\partial_{y_1}\phi-\partial_{y_2}k_2\partial_{y_2}\phi)(l_1(y_1)\partial_{y_1}v+l_2(y_1)\partial_{y_2}v) dy_1dy_2.
\end{eqnarray}
The coefficients in the quadratic term are slightly different from the ones in \eqref{basic}. However, similar argument still works in this case, which yields the following estimate
\begin{eqnarray}\no
\iint_{\Omega}(|\partial_{y_2}^2\phi|^2+|\partial_{y_1y_2}^2\phi|^2)dy_1 dy_2\leq \frac{C_*}{\sigma_*}\bigg(\iint_{\Omega}|\hat{F}|^2+(\partial_{y_2}\hat{F})^2dy_1 dy_2
+\displaystyle\sum_{j=2}^3\|g_j\|_1^2\bigg),
\end{eqnarray}
where the constant $C_*$ depends only on the $H^3$-norms of the coefficients $k_{ij}, k_i$ for $i,j=1,2$.

It follows from \eqref{LinearPotential3} that
\begin{eqnarray}\label{second-y}
\partial_{y_1}^2\phi=\hat{F}-k_{22}\partial_{y_2}^2\phi-2k_{12}\partial_{y_1y_2}^2 \phi-k_1\partial_{y_1}\phi-k_2\partial_{y_2}\phi,
\end{eqnarray}
which yields
\begin{eqnarray}\no
\iint_{\Omega}|\partial_{y_1}^2\phi|^2dy_1dy_2\leq \frac{C_*}{\sigma_*}\bigg(\iint_{\Omega}|\hat{F}|^2+(\partial_{y_2}\hat{F})^2dy_1dy_2+\|g_2\|_1^2+\|g_3\|_1^2\bigg).
\end{eqnarray}
Hence,
\begin{eqnarray}\label{second4}
\|\phi\|_{2}^2\leq\frac{C_*}{\sigma_*}\bigg(\|\hat{F}\|_{L^2}^2+\|\partial_{y_2}\hat{F}\|_{L^2}^2+\|g_2\|_1^2+\|g_3\|_1^2\bigg).
\end{eqnarray}

Rewrite the system \eqref{second1} as
\begin{eqnarray}\nonumber
\begin{cases}
\hat{L}v:=\displaystyle \sum_{i,j=1}^2 k_{ij} \partial_{y_i y_j}^2 v + k_1\partial_{y_1} v + k_2 \partial_{y_2} v =F_1,\\
r_0\partial_{y_1} v(r_0,y_2)+(r_0f'(r_0)-l_0)\partial_{y_2} v(r_0,y_2)=g'_2(y_2),\\
\partial_{y_2} v(r_1, y_2)=g_3'(y_2),
\end{cases}
\end{eqnarray}
where $F_1:=\partial_{y_2}\hat{F}-2 \partial_{y_2} k_{12} \partial_{y_1}v - \partial_{y_2} k_{22}\partial_{y_2}v-\partial_{y_2}k_1\partial_{y_1}\phi-\partial_{y_2}k_2\partial_{y_2}\phi$.

Applying the estimate \eqref{second4} to $v$ leads to
\begin{eqnarray}\nonumber
\|\partial_{y_2}\phi\|_{2}^2&\leq &\frac{C}{\sigma_*}(\|F_1\|_1^2+\|g_2\|_2^2+\|g_3\|_2^2)\\\nonumber
&\leq&\frac{C}{\sigma_*}\bigg(\|\hat{F}\|_2^2+\|\partial_{y_2}k_{22}\partial_{y_2}^2\phi\|_1^2+\|\partial_{y_2} k_{12}\partial_{y_1y_2}^2\phi\|_1^2+\|\partial_{y_2}k_1\partial_{y_1}\phi\|_1^2
\\\nonumber
&\quad&\quad+\|\partial_{y_2}k_2\partial_{y_2}\phi\|_1^2+\|g_2\|_2^2+\|g_3\|_2^2\bigg)\\\nonumber
&\leq&\frac{C}{\sigma_*}\bigg(\|\hat{F}\|_2^2+\|\partial_{y_2}k_{22}\|_2^2\|\partial_{y_2}^2\phi\|_1^2+\|\partial_{y_2} k_{12}\|_2^2\|\partial_{y_1y_2}^2\phi\|_1^2
+\|\partial_{y_2}k_1\|_2^2\|\partial_{y_1}\phi\|_1^2\\\nonumber
&\quad&\quad+\|\partial_{y_2}k_2\|_2^2\|\partial_{y_2}\phi\|_1^2+\|g_2\|_2^2+\|g_3\|_2^2\bigg),\\\nonumber
&\leq&\frac{C}{\sigma_*}\bigg(\|\hat{F}\|_2^2+\delta_0\|\partial_{y_2}\phi\|_2^2+\delta_0\|\phi\|_2^2+\|g_2\|_2^2+\|g_3\|_2^2\bigg),
\end{eqnarray}
which implies
\begin{align}\label{third1}
\|\partial_{y_2}\phi\|_{2}^2\leq\frac{C_*}{\sigma_*}(\|\hat{F}\|_2^2+\|g_2\|_2^2+\|g_3\|_2^2).
\end{align}

It follows from \eqref{second-y} that
\begin{eqnarray}\label{third-y}
\partial_{y_1}^3\phi&=&\partial_{y_1}\hat{F}-k_{22}\partial_{y_1}\partial_{y_2}^2\phi-2k_{12}\partial_{y_1}^2\partial_{y_2}\phi-k_1\partial_{y_1}^2\phi
-k_2\partial_{y_1y_2}^2\phi\\\nonumber
&\quad&-\partial_{y_1} k_{22}\partial_{y_2}^2\phi-2\partial_{y_1}k_{12}\partial_{y_1y_2}^2\phi-\partial_{y_1}k_1\partial_{y_1}\phi-\partial_{y_1}k_2\partial_{y_2}\phi,
\end{eqnarray}
which together with \eqref{second4} and \eqref{third1} yields
\begin{eqnarray}\label{third}
\|\phi\|_{3}^2\leq \frac{C_*}{\sigma_*}(\|\hat{F}\|_2^2+\|g_2\|_2^2+\|g_3\|_2^2).
\end{eqnarray}

Since $\p_{y_2}^2 \phi$ solves the following problem
\begin{align}\no
\begin{cases}
\hat{L}(\partial_{y_2}^2\phi):=\displaystyle \sum_{i,j=1}^2 k_{ij} \partial_{y_i y_j} (\partial_{y_2}^2\phi) + k_1\partial_{y_1} \partial_{y_2}^2\phi + k_2 \partial_{y_2} \partial_{y_2}^2\phi =F_2,\\
r_0\partial_{y_1}(\partial_{y_2}^2\phi)(r_0,y_2)+(r_0f'(r_0)-l_0)\partial_{y_2}(\partial_{y_2}^2\phi)(r_0,y_2)=g''_2(y_2),\\
\partial_{y_2}(\partial_{y_2}^2\phi)(r_1,y_2)=g''_3(y_2),
\end{cases}
\end{align}
where
\begin{eqnarray}\nonumber
F_2&=&\partial_{y_2}^2\hat{F}-\displaystyle \sum_{i,j=1}^2 \left(2\partial_{y_2} k_{ij} \partial_{y_i y_j}^2 \partial_{y_2}\phi+ \partial_{y_2}^2 k_{ij}\partial_{y_iy_j}\phi\right) -\sum_{i=1}^2 \left(2\partial_{y_2} k_i \partial_{y_iy_2}^2 \phi + \partial_{y_2}^2 k_i \partial_{y_i}\phi\right).
\end{eqnarray}

It follows from the derivation of \eqref{second4} that
\begin{align*}
\|\partial_{y_2}^2\phi\|_{2}^2\leq &\frac{C}{\sigma_*}(\|F_2\|_1^2+\|g_2\|_3^2+\|g_3\|_3^2)\\
\leq&\frac{C}{\sigma_*}\bigg(\|\hat{F}\|_3^2+\displaystyle\sum_{i,j=1}^2\left(\|\partial_{y_2}k_{ij}\partial_{y_i y_j}^2\partial_{y_2}\phi\|_1^2+\|\partial_{y_2}^2k_{ij}\partial_{y_i y_j}^2\phi\|_1^2\right)\\\nonumber
&+\sum_{i=1}^2\left( \|\partial_{y_2}k_i \partial_{y_iy_2}^2\phi\|_1^2+\|\partial_{y_2}^2k_i\partial_{y_i}\phi\|_1^2\right)+\|g_2\|_3^2+\|g_3\|_3^2\bigg)\\\nonumber
\leq&\frac{C}{\sigma_*}\bigg(\|\hat{F}\|_3^2+\|\partial_{y_2}^2 k_{22}\|_1^2\|\partial_{y_2}^2\phi\|_2^2+\|\partial_{y_2}^2k_{12}\|_1^2\|\partial_{y_1y_2}^2\phi\|_2^2
+\|\partial_{y_2}^2k_1\|_1^2\|\partial_{y_1}\phi\|_2^2\\
&+\|\partial_{y_2}^2k_2\|_1^2\|\partial_{y_2}\phi\|_2^2+\|\partial_{y_2}k_{22}\|_2^2\|\partial_{y_2}^3\phi\|_1^2+\|\partial_{y_2}k_{12}\|_2^2\|\partial_{y_2}^2\partial_{y_1}\phi\|_1^2
\\&+\|\partial_{y_2}k_1\|_2^2\|\partial_{y_1y_2}^2\phi\|_1^2+\|\partial_{y_2}k_2\|_2^2\|\partial_{y_2}^2\phi\|_1^2+\|g_2\|_3^2+\|g_3\|_3^2\bigg),\\
\leq&\frac{C_*}{\sigma_*}\bigg(\|\hat{F}\|_3^2+\delta_0\|\partial_{y_2}^2\phi\|_2^2+\delta_0\|\partial_{y_2}\partial_{y_1}^3\phi\|_0^2
+\delta_0\|\phi\|_3^2+\|g_2\|_3^2+\|g_3\|_3^2\bigg),
\end{align*}
where we have used the 2-D Sobolev embedding in the second and third inequalities.

This together with \eqref{third} implies
\begin{align}\no
\|\partial_{y_2}^2\phi\|_{2}^2\leq\frac{C_*}{\sigma_*}(\|\hat{F}\|_3^2+\delta_0\|\partial_{y_2}\partial_{y_1}^3\phi\|_0^2+\|g_2\|_3^2+\|g_3\|_3^2),
\end{align}
for suitably small $\delta_0$.

\eqref{third-y} implies that
\begin{eqnarray}\nonumber
&&\partial_{y_2}\partial_{y_1}^3\phi=\partial_{y_1y_2}^2\hat{F}-k_{22}\partial_{y_1}\partial_{y_2}^3\phi-2k_{12}\partial_{y_1}^2\partial_{y_2}^2\phi
-\partial_{y_1} k_{22}\partial_{y_2}^3\phi-2\partial_{y_1}k_{12}\partial_{y_1}\partial_{y_2}^2\phi\\\nonumber
&&\quad\quad- \partial_{y_2}k_{22}\partial_{y_1}\partial_{y_2}^2\phi-2\partial_{y_2}k_{12}\partial_{y_1}^2\partial_{y_2}\phi
-\partial_{y_1y_2}^2 k_{22}\partial_{y_2}^2\phi-2\partial_{y_1y_2}^2k_{12}\partial_{y_1y_2}^2\phi
\\\nonumber
&&\quad\quad-\sum_{i=1}^2 \left(\partial_{y_2} k_i \partial_{y_i y_1}^2 \phi + k_i\partial_{y_i}\partial_{y_1y_2}^2\phi\right)-\sum_{i=1}^2\left(\partial_{y_1}k_i\partial_{y_iy_2}\phi+\partial_{y_1y_2}^2k_i\partial_{y_i}\phi\right),
\end{eqnarray}
from which one can derive that
\begin{eqnarray}\no
\|\partial_{y_2}\partial_{y_1}^3\phi\|_{L^2}^2\leq\frac{C_*}{\sigma_*}(\|\hat{F}\|_3^2+\delta_0\|\partial_{y_2}\partial_{y_1}^3\phi\|_0^2+\|g_2\|_3^2+\|g_3\|_3^2).
\end{eqnarray}
Thus
\begin{align}\no
\|\partial_{y_2}\partial_{y_1}^3\phi\|_{L^2}^2\leq\frac{C_*}{\sigma_*}(\|\hat{F}\|_3^2+\|g_2\|_3^2+\|g_3\|_3^2)
\end{align}
for suitably small $\delta_0$.

It remains to estimate $\partial_{y_1}^4\phi$. Since $\partial_{y_1}^4\phi=\partial_{y_1}^2(\hat{F}-k_{22}\partial_{y_2}^2\phi-2k_{12}\partial_{y_1y_2}^2\phi-k_1\partial_{y_1}\phi-k_2\partial_{y_2}\phi)$, it holds that
\begin{align}\no
\|\partial_{y_1}^4\phi\|_{L^2}^2\leq\frac{C_*}{\sigma_*}(\|\hat{F}\|_3^2+\|g_2\|_3^2+\|g_3\|_3^2).
\end{align}
In summary, we obtain
\begin{align}\no
\|\phi\|_{4}^2\leq \frac{C_*}{\sigma_*}(\|\hat{F}\|_3^2+\|g_2\|_3^2+\|g_3\|_3^2).
\end{align}

\end{proof}

\subsubsection{Proof of Theorem \ref{main1}} \label{existence}\noindent

We then turn to the proof of Theorem \ref{main1}. Inspired by Kuzmin\cite{Kuzmin2002}, we use Galerkin's method to construct approximate solutions and a simple contraction mapping argument yields the solution to the nonlinear problem. Slightly different from the case in \cite{Kuzmin2002}, the coefficient of $\p_{y_2}^2 \phi$ in \eqref{LinearPotential3} changes sign, not that of $\p_{y_1}^2 \phi$, it is not necessary to add a third order dissipative term to \eqref{LinearPotential3} to construct the approximate solutions.

Since the solution is periodic in $y_2$, it is natural to use the Fourier series to construct the approximate solution to the linearized problem. Note that the $H^4(\Omega)$ energy estimate is obtained only for the linearized problem with smooth coefficients, one needs to mollify the coefficients in \eqref{LinearPotential3}.

Since $k_{12}$, $k_{22}$, $k_1$, $k_2$, $\hat{F}\in H^3$, there exists sequences of smooth functions $\{k_{12}^\eta\}$, $\{k_{22}^\eta\}$, $\{k_1^\eta\}$, $\{k_2^\eta\}$, $\{\hat{F}^\eta\}$ such that $\|k_{12}^\eta-k_{12}\|_3\rightarrow0$, $\|k_{22}^\eta-k_{22}\|_3\rightarrow0$, $\|k_1^\eta-k_1\|_3\rightarrow0$, $\|k_2^\eta-k_2\|_3\rightarrow0$, $\|\hat{F}^\eta-\hat{F}\|_3\rightarrow0$ as $\eta\rightarrow0$. Consider the linear boundary value problem
\begin{eqnarray}\label{approx1}\begin{cases}
\hat{L}^\eta\phi\equiv \displaystyle \sum_{i,j=1}^2 k_{ij}^{\eta} \partial_{y_i y_j}^2\phi+k_1^\eta\partial_{y_1}\phi+k_2^\eta\partial_{y_2}\phi= \hat{F}^{\eta},\\
r_0\partial_{y_1}\phi(r_0,y_2)+(r_0f'(r_0)-l_0)\partial_{y_2}\phi(r_0,y_2)=g_2(y_2),\\
\partial_{y_2}\phi(r_1,y_2)=g_3(y_2), \\
\phi(r_1, f(r_1))=0.
\end{cases}\end{eqnarray}

Note that $\int_0^{2\pi} g_3(y_2) dy_2=0$, so one may assume that $g_3\equiv 0$, otherwise consider the function $\phi-\int_0^{y_2} g_3(t)dt$. Thus the boundary condition on $y_1=r_1$ becomes $\phi(r_1,y_2)=0$ for any $y_2\in \mathbb{T}_{2\pi}$. Choose the standard orthonormal basis $\{h_j(y_2)\}_{j=1}^{\infty}$ of $L^2(\mathbb{T}_{2\pi})$, where for each positive integer $m\in\mathbb{N}$:
\be\no
h_1(y_2)=\frac{1}{\sqrt{2\pi}},\ \ h_{2m}(y_2)=\frac{1}{\sqrt{\pi}}\sin (m y_2),\ \ h_{2m+1}(y_2)= \frac{1}{\sqrt{\pi}}\cos (m y_2),\cdots.
\ee
and we construct approximate solutions to \eqref{approx1} of the form
\begin{align}\no
\phi^{N,\eta}(y_1,y_2)=\sum_{j=1}^{2N+1} A_j^{N,\eta}(y_1)h_j(y_2),
\end{align}
where $A_j^{N,\eta}(y_1)$ are determined by the system of $2N+1$ second-order ordinary differential equations supplemented with $2(2N+1)$ boundary conditions:
\begin{align}\label{ODE1}
\begin{cases}
\int_0^{2\pi}(\hat{L}^\eta\phi^{N,\eta}-\hat{F}^\eta)h_m(y_2)dy_2=0,\ \ \ m=1,2,...,2N+1,\\
\int_0^{2\pi}(r_0\partial_{y_1}\phi^{N,\eta}(r_0,y_2)+(r_0f'(r_0)-l_0)\partial_{y_2}\phi^{N,\eta}(r_0,y_2))h_m(y_2)dy_2=\int_0^{2\pi}g_2(y_2) h_m(y_2) dy_2,\\
\int_0^{2\pi} \phi^{N,\eta}(r_1,y_2))h_m(y_2)dy_2=0,\ \ \ m=1,2,...,2N+1.
\end{cases}
\end{align}

Define
\be\no
&&a_{jm}^{N,\eta}= \int_0^{2\pi} (2k_{12}^{\eta}(y_1,y_2) h_j'(y_2)+ k_1^{\eta}(y_1,y_2) h_j(y_2)) h_m(y_2) dy_2,\\\no
&&b_{jm}^{N,\eta}= \int_0^{2\pi} (k_{22}^{\eta} h_j''(y_2)+k_2^{\eta} h_j'(y_2)) h_m(y_2) dy_2,\ c_{jm}=\int_0^{2\pi} h_j'(y_2) h_m(y_2) dy_2.
\ee
Then the system \eqref{ODE1} reduces to
\be\label{ODE2}\begin{cases}
\displaystyle\frac{d^2}{dy_1^2} A_m^{N,\eta}(y_1) + \sum_{j=1}^{2N+1} a_{jm}^{N,\eta} \frac{d}{dy_1} A_j^{N,\eta}(y_1) + \sum_{j=1}^{2N+1} b_{jm}^{N,\eta} A_j^{N,\eta}(y_1)=\int_0^{2\pi} \hat{F}^{\eta} h_m(y_2) dy_2,\\
\displaystyle r_0\frac{d}{dy_1} A_m^{N,\eta}(r_0)+(r_0f'(r_0)-l_0) \sum_{j=1}^{2N+1} c_{jm} A_j^{N,\eta}(r_0)= \int_0^{2\pi} g_2(y_2) h_m(y_2) dy_2,\\
A_m^{N,\eta}(r_1)=0,\ \ m=1,\cdots, 2N+1.
\end{cases}\ee

Using the functions $l_1(y_1), l_2(y_1)$ defined in Section \ref{energy}, multiplying the first equation, the $(2m)^{th}$ equation and the $(2m+1)^{th}$ equation in \eqref{ODE2} by $l_1(y_1)\frac{d}{dy_1} A_1^{N,\eta}$,  $l_1(y_1)\frac{d}{dy_1} A_{2m}^{N,\eta}- m l_2(y_1)A_{2m+1}^{N,\eta}$ and $l_1(y_1)\frac{d}{dy_1} A_{2m+1}^{N,\eta}+ m l_2(y_1)A_{2m}^{N,\eta}$ respectively, summing from $1$ to $2N+1$, one can get after integrating over $[r_0,r_1]$ that
\begin{eqnarray}\label{ODE4}
\iint_{\Omega}(\hat{L}^\eta\phi^{N,\eta}-\hat{F}^\eta)(l_1(y_1)\partial_{y_1} \phi^{N,\eta}+ l_2(y_1)\partial_{y_2} \phi^{N,\eta})dy_1 dy_2=0.
\end{eqnarray}
Similar argument as in Lemma \ref{basic-energy} yields
\begin{align}\label{approx3}
\|\phi^{N,\eta}\|_{1}\leq \frac{C_*}{\sigma_*}(\|\hat{F}^{\eta}\|_{L^2}+\sum_{j=2,3}\|g_j\|_{L^2(\mathbb{T}_{2\pi})}).
\end{align}
The following higher order derivatives estimate can be derived similarly as in Lemma \ref{highorder}:
\begin{align}\label{approx2}
\|\phi^{N,\eta}\|_{4}\leq \frac{C_*}{\sigma_*}(\|\hat{F}^{\eta}\|_{3}+\|g_2\|_3+ \|g_3\|_3).
\end{align}

This estimate implies the uniqueness of the solution to Problem \eqref{ODE2}. By Fredholm alternative theorem for second order elliptic systems, the uniqueness ensures the existence of the solution to \eqref{ODE2}. Since the coefficients of the \eqref{ODE2} are smooth, so the solutions $\phi^{N,\eta}$ are smooth. It should be emphasized that the estimates obtained in Section \ref{energy} only involves the $H^3$ norm of the coefficients, therefore the bound on the right hand side of \eqref{approx2} is uniformly in $N, \eta$. For any fixed $\eta>0$, by the weak compactness of a bounded set in $H^4(\Omega)$, there exists a subsequence $\{\phi^{N_j,\eta}\}_{j=1}^{\infty}$ that converges weakly to $\phi^{\eta}$ in $H^4$ and the convergence is strong in $H^3(\Omega)$. Therefore $\phi^{\eta}$ will be the unique solution to \eqref{approx1}. The estimate \eqref{approx2} also holds for $\phi^{\eta}$ with a constant $C_*$ independent of $\eta$, from which we can find a subsequence $\{\phi^{\eta_j}\}_{j=1}^{\infty}$ converging weakly to a function $\phi$ in $H^4(\Omega)$. In conclusion, we have proved that the problem \eqref{LinearPotential3} has a unique solution $\phi\in H^4$ with the estimate
\begin{align}\label{energyest1}
\|\phi\|_{4}\leq \frac{C_*}{\sigma_*}(\|\hat{F}\|_3+\|g_2\|_3+\|g_3\|_3)\leq C(\epsilon+\delta_0^2).
\end{align}

Hence the mapping $\mathcal{T}$ is well-defined in $\mathcal{X}$ for sufficiently small $\delta_0=\sqrt{\eps}$. It remains to show that the mapping $\mathcal{T}$ is contractive in a low order norm for sufficiently small $\delta_0$. Suppose that $\phi^{(i)}=\mathcal{T}\bar{\phi}^{(i)} (i=1,2)$ for any $\bar{\phi}^{(1)}$, $\bar{\phi}^{(2)}\in \mathcal{X}$. Then for $k=1,2$,
\begin{align}\no
\begin{cases}
\hat{L}^{(k)}\phi^{(k)}\equiv \displaystyle \sum_{i,j=1}^2 k_{ij}(\bar{U}_1^{(k)},\bar{U}_2^{(k)})\partial_{y_i y_j}^2\phi^{(k)}+\sum_{i=1}^2 k_i(\bar{U}_1^{(k)},\bar{U}_2^{(k)})\partial_{y_i}\phi^{(k)}=\hat{F}(\bar{U}_1^{(k)},\bar{U}_2^{(k)}),\\
r_0\partial_{y_1}\phi^{(k)}(r_0, y_2)+(r_0f'(r_0)-l_0)\partial_{y_2}\phi^{(k)}(r_0,y_2)=g_2(y_2),\\
\partial_{y_2}\phi^{(k)}(r_1,y_2)=g_3(y_2),\\
\phi^{(k)}(r_1,f(r_1))=0.
\end{cases}
\end{align}
Thus
\begin{align}\no
\begin{cases}
\hat{L}^{(1)}(\phi^{(1)}-\phi^{(2)})
=\hat{F}(\bar{U}_1^{(1)},\bar{U}_2^{(1)})-\hat{F}(\bar{U}_1^{(2)},\bar{U}_2^{(2)})-(\hat{L}^{(1)}-\hat{L}^{(2)})\phi^{(2)}\\
r_0(\partial_{y_1}\phi^{(1)}-\partial_{y_1}\phi^{(2)})(r_0,y_2)+(r_0f'(r_0)-l_0)(\partial_{y_2}\phi^{(1)}-\partial_{y_2}\phi^{(2)})(r_0,y_2)=0,\\
(\phi^{(1)}-\phi^{(2)})(r_1,y_2)=0,
\end{cases}
\end{align}
which implies
\begin{align}\no
\|\mathcal{T}\bar{\phi}^{(1)}-\mathcal{T}\bar{\phi}^{(2)}\|_1=\|\phi^{(1)}-\phi^{(2)}\|_1&\leq C_*\|\hat{F}(\bar{U}_1^{(1)},\bar{U}_2^{(1)})-\hat{F}(\bar{U}_1^{(2)},\bar{U}_2^{(2)})-(\hat{L}^{(1)}-\hat{L}^{(2)})\phi^{(2)}\|_0\\\no
&\leq C_* \delta_0\|\bar{\phi}^{(1)}-\bar{\phi}^{(2)}\|_1\leq \frac{1}{2}\|\bar{\phi}^{(1)}-\bar{\phi}^{(2)}\|_1,
\end{align}
since $\mathcal{T}\bar{\phi}^{(2)}$, $\bar{\phi}^{(1)}$, $\bar{\phi}^{(2)}\in \mathcal{X}$, and $\hat{F}(U_1,U_2)$, $k_{12}(U_1,U_2), k_{22}(U_1,U_2)$, $k_i(U_1,U_2)$, $i=1,2$ are smooth functions of $(U_1,U_2)$, and one can choose $\delta_0$ small enough such that $\mathcal{T}$ is a contractive mapping in $H^1$-norm. Then there exists a unique $\phi\in \mathcal{X}$ to $\mathcal{T}\phi=\phi$.


In conclusion, we have shown that there exists a small $\epsilon_0>0$ such that for any $0<\epsilon<\epsilon_0$, the problem \eqref{Pro} has a unique solution in $H^4(\Omega)$ with the estimate $\|\phi\|_4\leq C_*\epsilon$. That is, the background transonic flow is structurally stable within irrotational flows under perturbations of the flow angles at the inner and outer circular cylinder.

Finally, we examine the location of all the sonic points which satisfy $|{\bf M}(r,\th)|^2=1$, where ${\bf M}=(M_1,M_2)^t:=(\frac{U_1}{c(\rho)},\frac{U_2}{c(\rho)})^t$. It follows from \eqref{estimate1} and the Sobolev embedding $H^3(\Omega)\hookrightarrow C^{1,\alpha}(\overline{\Omega})$ for any $\alpha\in (0,1)$ that
\be\no
\||{\bf M}|^2-|{\bf M}_b|^2\|_{C^{1,\alpha}(\overline{\Omega})}\leq \||{\bf M}|^2-|{\bf M_b}|^2\|_3\leq C_*\epsilon.
\ee
Note that
\be\no
|{\bf M}_b(r_0)|^2>1, \ |{\bf M}_b(r_1)|^2<1, \ \ \displaystyle \sup_{r\in[r_0,r_1]} \frac{d}{d r}|{\bf M}_b(r)|^2<0.
\ee
Thus for sufficiently small $\epsilon$,  $|{\bf M}(r_0,\theta)|^2>1, \ |{\bf M}(r_1,\theta)|^2<1$ for any $\theta\in\mathbb{T}_{2\pi}$ and $\frac{\p}{\p r}|{\bf M}(r,\theta)|^2<0$ for any $(r,\theta)\in\Omega$. Therefore for each $\theta\in \mathbb{T}_{2\pi}$, there exists a unique $s(\th)\in (r_0,r_1)$ such that $|{\bf M}(s(\th),\th)|^2=1$. Also by the implicit function theorem, the function $s\in C^{1}(\mathbb{T}_{2\pi})$. Furthermore, since
\be\no
||{\bf M}_b(s(\theta))|^2-|{\bf M}_b(r_c)|^2|&=&|{\bf M}_b(s(\theta))|^2-|{\bf M}(s(\theta),\theta)|^2|\\\no
&\leq& \||{\bf M}|^2-|{\bf M}_b|^2\|_{C^{1,\alpha}(\overline{\Omega})} \leq C_*\epsilon,
\ee
one can deduce that $|s(\theta)-r_c|\leq C_*\epsilon$ for any $\theta\in \mathbb{T}_{2\pi}$. Differentiating the identity $|{\bf M}(s(\th),\th)|^2=1$ with respect to $\theta$ yields
\be\no
s'(\theta)=-\left(\frac{\p}{\p r}|{\bf M}|^2 (s(\theta),\theta) \right)^{-1}\frac{\p}{\p\theta}|{\bf M}|^2 (s(\theta),\theta)
\ee
and the estimate \eqref{estimate2} holds.

\subsection{Proof of Theorem \ref{2dmain}}\label{proofTh12}\noindent

We now turn to the general case that the flow may be rotational and prove Theorem \ref{2dmain}. It is well-known that the steady Euler system is coupled elliptic-hyperbolic in subsonic region and changes type when the flow changes from subsonic to supersonic. There are several different decompositions developed by different authors for different purposes \cite{ccs06,ccf07,chen08,lxy09b,lxy13,xy05}. We will employ the deformation-curl decomposition developed in \cite{WengXin19,weng2019} to deal with the elliptic-hyperbolic coupled structure in the steady Euler equations. The Bernoulli's law yields
\be\label{2d-density}
\rho=H(|{\bf U}|^2,B,A)=\left(\frac{\gamma-1}{A\gamma} (B-\frac{1}{2} |{\bf U}|^2)\right)^{\frac{1}{\gamma-1}}.
\ee
As discussed in \cite{WengXin19,weng2019}, one can show that if a smooth flow does not contain the vacuum and the stagnation points, the steady Euler system \eqref{2d} is equivalent to the following system
\begin{eqnarray}\label{2d1}
\begin{cases}
(c^2-U_1^2)\p_r U_1+(c^2-U_2^2)\frac1r\p_\th U_2-U_1U_2(\p_rU_2+\fr1r\p_\th U_1)+\fr1r c^2 U_1\\
=-(U_1\p_r+U_2\fr1r\p_\th)B+\fr{c^2}{(\ga-1)A}(U_1\p_r+U_2\fr1r\p_\th)A,\\
\frac{1}{r}U_2(\p_\th U_1-\p_r(rU_2))=\frac{ (B-\frac{1}{2}|{\bf U}|^2)}{A\gamma}\p_r A-\p_r B,\\
(U_1\partial_r +\frac{U_2}{r}\partial_{\theta}) B=0,\\
(U_1\partial_r +\frac{U_2}{r}\partial_{\theta}) A=0,
\end{cases}
\end{eqnarray}
where the first equation in \eqref{2d1} is derived by substituting \eqref{2d-density} into the density equation.

Define $\hat{U}_1$ and $\hat{U}_2$ as in \eqref{difference} and set $\hat{B}= B-B_0$ and $\hat{A}= A-A_0$. Then ${\bf \hat{U}}$, $\hat{B}$ and $\hat{A}$ satisfy
\begin{eqnarray}\label{aa}
\begin{cases}
A_{11}\p_r \hat{U}_1+rA_{22}\p_\th \hat{U}_2+r A_{12}\p_r\hat{U}_2+ A_{21}\p_\theta \hat{U}_1+e_1(r)\hat{U}_1+\tilde{e}_2(r)\hat{U}_2=F_1({\bf U},B,A),\\
\frac{1}{r}\p_\th \hat{U}_1-\frac{1}{r}\p_r(r\hat{U}_2)=F_2({\bf U},B,A),\\
(U_1\partial_r +\frac{U_2}{r}\partial_{\theta}) \hat{B}=0,\\
(U_1\partial_r +\frac{U_2}{r}\partial_{\theta})\hat{A}=0,
\end{cases}
\end{eqnarray}
where
\begin{eqnarray}\label{2dboundary-coe}
\begin{cases}
A_{11}({\bf U},B)=c^2(B, |{\bf U}|^2)-U_{1}^2,\ \ A_{22}({\bf U},B)=\frac{1}{r^2}(c^2(B, |{\bf U}|^2)-U_{2}^2),\\
A_{12}=A_{21}=-\frac{U_{1}U_{2}}{r},\ \ \ c^2(B, |{\bf U}|^2)= (\gamma-1)(B-\frac{1}{2} |{\bf U}|^2),\\
e_1(r)=\frac{c^2(\rho_b)+U_{b2}^2}r +(\ga+1)\frac{1+M_{b2}^2}{r(1-M_{b1}^2)}U_{b1}^2-\frac{(\ga-1)U_{b1}^2}r,\\
\tilde{e}_2(r)=\frac{(\ga-1)(1+M_{b2}^2)}{r(1-M_{b1}^2)}U_{b1}U_{b2}-\frac{\ga-2}rU_{b1}U_{b2}=\frac1r U_{b1}U_{b2}(1+(\ga-1)\frac{|{\bf M}_b|^2}{1-M_{b1}^2}),\\
F_1({\bf U},B,A)=-\frac1r (c^2-c_b^2) \hat{U}_1-((\ga-1)(\hat{B}-\fr12\hat{U}_2^2)-\frac{\ga+1}{2}\hat{U}_1^2)U_{b1}'(r)\\\quad\quad\quad\quad\quad\quad\ +\hat{U}_1\hat{U}_2U_{b2}'(r)
-\frac1r (\ga-1)(\hat{B}-\frac12\hat{U}_1^2-\frac12\hat{U}_2^2) U_{b1}\\\quad\quad\quad\quad\quad\quad\ -(U_1\p_r+U_2\fr1r\p_\th)\hat{B}+\fr{c^2}{(\ga-1)A}(U_1\p_r+U_2\fr1r\p_\th)\hat{A},\\
F_2({\bf U},B,A)=\frac{1}{U_2}(\frac{(B-\frac{1}{2}|{\bf U}|^2)}{A\gamma}\p_r\hat{A}-\p_r \hat{B}).
\end{cases}
\end{eqnarray}
Note that $e_1(r)$ is same as the one defined in Section \ref{key} and $\tilde{e}_2(r)=r e_2(r)-\frac{U_{b1}U_{b2}}{r}$.

The boundary conditions in \eqref{2dboundary1}-\eqref{2dboundary4} reduce to
\be\no\begin{cases}
\hat{B}(r_1,\th)=\eps B_1(\th),\ \ \ \hat{A}(r_1,\theta)=\epsilon A_1(\theta),\\
\hat{U}_1(r_0,\th)-l_0 \hat{U}_2(r_0,\theta)=\epsilon g_{0}(\theta),\\
\hat{U}_2(r_1,\th)=\epsilon g_{1}(\theta).
\end{cases}\ee

\br\label{loss}
{\it Note that the first two equations in \eqref{aa} will be regarded as first order mixed type equations for $U_1$ and $U_2$, the energy estimates in previous section indicate that the regularity of the solutions $U_1, U_2$ would be at best the same as the source terms on the right hand side. Hence if one looks for the solution $({\bf U}, B, A)$ in $(H^3(\Omega))^4$, then $F_2({\bf U},B,A)$ belongs only to $H^2(\Omega)$ and there appears a loss of derivatives. Similar issue occurs for the structural stability of 1-D supersonic flows to the steady Euler-Poisson system \cite{bdxx19}. However, the approach in \cite{bdxx19} can not be adapted directly here, since one can not prescribe the boundary data for all the flow quantities at the entrance in this case. To overcome this difficulty, besides introducing the stream function to solve the transport equations, we also observe that the regularity of the flows in the subsonic region can be improved if the data at the entrance have better regularity. Thus we will choose some appropriate functional spaces and design an elaborate two-layer iteration scheme to prove Theorem \ref{2dmain}.
}\er

Define $\Omega_{ue}=\{(r,\theta): \frac12 r_c+\frac12 r_1< r< r_1, \theta\in [0,2\pi]\}$ and
\be\no
&&\Omega_{lue}=\{(r,\theta): \fr34 r_c+\fr14 r_1< r< r_1, \theta\in [0,2\pi]\},\\\no
&&\tilde{\Omega}_{lue}=\{(r,\theta): \fr78 r_c+\fr18 r_1< r< r_1, \theta\in [0,2\pi]\}.
\ee
Then $\Omega_{ue}\subset \Omega_{lue}\subset \tilde{\Omega}_{lue}$. Set
\begin{align}\nonumber
\mathcal{X}_1=&\{{\bf U}(r,\th)\in (H^3(\Omega))^2: \|\hat{{\bf U}}\|_{H^3(\Omega)}\leq \delta_0\},\\\no
\mathcal{X}_2=&\{(B(r,\th),A(r,\th)) \in (H^4(\Omega))^2\cap(C^{3,\alpha}(\overline{\Omega_{lue}}))^2\cap(C^{4,\alpha}(\overline{\Omega_{ue}}))^2:\\\no
& \quad\quad\quad\|(\hat{B},\hat{A})\|_{H^4(\Omega)}+ \|(\hat{B},\hat{A})\|_{C^{3,\alpha}(\overline{\Omega_{lue}})}+ \|(\hat{B},\hat{A})\|_{C^{4,\alpha}(\overline{\Omega_{ue}})}\leq \delta_1\}
\end{align}
with positive constants $\delta_0,\delta_1>0$ to be specified later. For fixed $(\bar{B},\bar{A})\in\mathcal{X}_2$ and for any function $\bar{{\bf U}}\in \mathcal{X}_1$, we first construct an operator $\mathcal{T}^{(\bar{B},\bar{A})}$: $\bar{{\bf U}}\in \mathcal{X}_1\mapsto{\bf U}\in \mathcal{X}_1$, where ${\bf U}={\bf U}_b+ \hat{{\bf U}}$ is obtained by resolving the following boundary value problem
\begin{eqnarray}\label{7step2}
\begin{cases}
A_{11}(\bar{{\bf U}},\bar{B},\bar{A})\p_r \hat{U}_1+r A_{22}(\bar{{\bf U}},\bar{B},\bar{A})\p_\th \hat{U}_2+r A_{12}(\bar{{\bf U}})\p_r\hat{U}_2+ A_{21}(\bar{{\bf U}})\p_\th \hat{U}_1\\
\quad\quad+e_1(r)\hat{U}_1+\tilde{e}_2(r)\hat{U}_2 =F_1(\bar{{\bf U}},\bar{B},\bar{A}),\\
\frac1r\p_\th \hat{U}_1-\frac{1}{r}\p_r(r\hat{U}_2)=F_2(\bar{{\bf U}},\bar{B},\bar{A}),\\
\hat{U}_1(r_0,\th)-l_0 \hat{U}_2(r_0,\theta)=\epsilon g_{0}(\theta),\\
\hat{U}_2(r_1,\th)=\epsilon g_{1}(\theta).
\end{cases}
\end{eqnarray}
Note that the equations in \eqref{7step2} form a linear first order mixed type system with coefficients given in \eqref{2dboundary-coe}.

Since $\bar{B}, \bar{A}\in \mathcal{X}_2, \bar{{\bf U}}\in \mathcal{X}_1$, there holds that
\begin{eqnarray}\no
\|F_1(\bar{{\bf U}},\bar{B},\bar{A})\|_{H^3(\Omega)}\leq C_0(\delta_1+\delta_0^2),\ \ \ \|F_2(\bar{{\bf U}},\bar{B},\bar{A})\|_{H^3(\Omega)}\leq C_0\delta_1.
\end{eqnarray}

Let $\phi_1(r,\theta)$ be the unique solution to the following problem
\begin{align}\no
\begin{cases}
(\p_r^2+\frac{1}{r}\p_r+ \frac{1}{r^2}\p_\th^2)\phi_1=F_2(\bar{{\bf U}},\bar{B},\bar{A})\in H^3(\Omega),\\
\phi_1(r_0,\theta)=\phi_1(r_1,\theta)=0.
\end{cases}
\end{align}
Then $\phi_1(r,\th)\in H^{5}(\Omega)$ and satisfies the estimate
\begin{align}\no
\|\phi_1\|_{H^{5}(\Omega)}\leq C_0\|F_2(\bar{{\bf U}},\bar{B},\bar{A})\|_{H^3(\Omega)}\leq C_0\delta_1.
\end{align}

Define $V_1= \hat{U}_1-\frac{1}{r}\p_{\theta}\phi_1$ and $V_2= \hat{U}_2+\p_{r}\phi_1$. Then
\begin{eqnarray}\label{7step21}
\begin{cases}
A_{11}(\bar{{\bf U}},\bar{B},\bar{A})\p_r V_1+rA_{22}(\bar{{\bf U}},\bar{B},\bar{A})\p_\th V_2+r A_{12}(\bar{{\bf U}})\p_rV_2+A_{21}(\bar{{\bf U}})\p_\th V_1\\
+e_1(r)V_1+\tilde{e}_2(r)V_2=F_3(\bar{U}_1,\bar{U}_2,\bar{B},\bar{A}),\\
\frac{1}{r}\p_\theta V_1-\frac{1}{r}\p_r(r V_2)=0,\\
V_1(r_0,\th)-l_0 V_2(r_0,\theta)=\epsilon g_{0}(\theta)-\frac{1}{r}\p_{\theta}\phi_1(r_0,\theta)- l_0\p_{r} \phi_1(r_0,\theta),\\
V_2(r_1,\th)=\epsilon g_{1}(\theta)+ \p_r\phi_1(r_1,\theta),
\end{cases}
\end{eqnarray}
where
\be\no
&&F_3(\bar{{\bf U}},\bar{B},\bar{A})= F_1(\bar{{\bf U}},\bar{B},\bar{A})- A_{11}(\bar{{\bf U}},\bar{B},\bar{A}) \p_r(\frac{1}{r}\p_{\theta}\phi_1) +r A_{22}(\bar{{\bf U}},\bar{B},\bar{A})\p_{r\theta}^2\phi_1 \\\no
&&\quad+ r A_{12}(\bar{\bf U})\p_{r}^2\phi_1 -\frac{A_{21}(\bar{\bf U})}{r} \p_{\theta}^2\phi_1- \frac{e_1(r)}{r} \p_{\theta}\phi_1 + \tilde{e}_2(r)\p_r\phi_1\in H^3(\Omega),
\ee
and
\begin{align}\no
\|F_3(\bar{{\bf U}},\bar{B},\bar{A})\|_{H^3(\Omega)}\leq C_0(\delta_1+\delta_0^2).
\end{align}
Introduce the potential function
\be\no
\phi_2(r,\theta)= \int_{r_1}^r V_1(\tau,\theta)d\tau + \int_0^{\theta} r_1 (V_2(r_1,s)+\tilde{d}_0) d\tau,
\ee
where $\tilde{d}_0= -\frac{1}{2\pi}\int_0^{2\pi} [\epsilon g_1(\theta)+ \p_r \phi_1(r_1,\theta)]d\theta$ is introduced to guarantee that $\phi_2(r,\theta)=\phi_2(r,\theta+2\pi)$. Then $\phi_2$ is periodic in $\theta$ with period $2\pi$ and satisfies
\be\label{phi2}
\p_r \phi_2(r,\theta) = V_1(r,\theta),\ \ \ \p_{\theta} \phi_2(r,\theta) =r V_2(r,\theta)+ r_1 \tilde{d}_0.
\ee

Substituting \eqref{phi2} into \eqref{7step21} leads to the following boundary value problem for a second-order linear mixed type equation
\begin{eqnarray}\label{7step22}
\begin{cases}
A_{11}(\bar{{\bf U}},\bar{B},\bar{A})\p_r^2\phi_2+A_{22}(\bar{{\bf U}},\bar{B},\bar{A})\p_\th^2\phi_2+(A_{12}+A_{21})(\bar{{\bf U}})\p^2_{r\theta}\phi_2\\\quad\quad\quad+e_1(r)\p_r\phi_2+E_2(\bar{{\bf U}})\p_\th\phi_2
=F_4(\bar{{\bf U}},\bar{B},\bar{A}),\\
\p_r\phi_2(r_0,\theta)-\frac{l_0}{r_0}\p_{\theta}\phi_2(r_0,\theta)=\tilde{g}_0(\theta),\\
\p_\th\phi_2(r_1,\th)=\tilde{g}_1(\theta),\ \ \ \phi_2(r_1,0)=0,
\end{cases}
\end{eqnarray}
where
\be\no
&&F_4(\bar{{\bf U}},\bar{B},\bar{A})=F_3(\bar{{\bf U}},\bar{B},\bar{A})+ \left(-\frac{A_{12}}{r}+ \frac{\tilde{e}_2(r)}{r}\right) r_1 \tilde{d}_0\in H^3,\\\no
&&E_2(\bar{{\bf U}})= \frac{\tilde{e}_2(r)}{r} -\frac{A_{12}(\bar{{\bf U}})}{r}=\frac{\tilde{e}_2(r)}{r}+ \frac{\bar{U}_1\bar{U}_2}{r^2},\\\no
&&\tilde{g}_0(\theta)=\epsilon g_{0}(\theta)-\frac{1}{r}\p_{\theta}\phi_1(r_0,\theta)- l_0\p_{r} \phi_1(r_0,\theta)-\frac{l_0 r_1}{r_0} \tilde{d}_0,\\\no
&&\tilde{g}_1(\theta)=\epsilon r_1 g_1(\theta)+r_1\p_r\phi_1(r_1,\th)+ r_1\tilde{d}_0
\ee
and
\be\no
\|E_2(\bar{{\bf U}})-e_2(r)\|_{H^3(\Omega)}=\|\frac{1}{r^2}(\bar{U}_1\bar{U}_2-U_{b1}U_{b2})\|_{H^3(\Omega)}\leq C_0\delta_0.
\ee

The problem \eqref{7step22} is only different slightly from \eqref{Pro}, one can adapt the same ideas in previous section to show the existence and uniqueness of a smooth solution $\phi_2\in H^4(\Omega)$ to \eqref{7step22} with the estimate
\be\label{7step23}
\|\phi_2\|_{H^4(\Omega)}\leq C(\|F_4(\bar{{\bf U}},\bar{B},\bar{A})\|_{H^3(\Omega)}+ \sum_{j=0,1}\|\tilde{g}_j\|_{H^3(\mathbb{T}_{2\pi})})\leq C_0(\epsilon+\delta_1+\delta_0^2).
\ee
Therefore, $\mathcal{T}^{(\bar{B},\bar{A})}$ is well-defined if one sets $\delta_0=\sqrt{\eps+\delta_1}$ and selects $\sqrt{\eps+\delta_1}\leq \frac{1}{2C_0}$.

It remains to show that the iteration mapping $\mathcal{T}^{(\bar{B},\bar{A})}$ is contractive in a low order norm for sufficiently small $\eps+\delta_1$. Set ${\bf U}^{(i)}=\mathcal{T}^{(\bar{B},\bar{A})}(\bar{{\bf U}}^{(i)}) (i=1,2)$ for any $\bar{\U}^{(1)}$, $\bar{\U}^{(2)}\in \mathcal{X}_1$ and denote $\bar{U}_j^{(1)}-\bar{U}_j^{(2)}$ by $\bar{V}_j$, $j=1,2$, $U_j^{(1)}-U_j^{(2)}$ by $V_j$, $j=1,2$. Then, it follows from \eqref{7step2} that
\begin{eqnarray}\no
\begin{cases}
A_{11}^{(1)}\p_rV_1+r A_{22}^{(1)}\p_{\theta} V_2+ r A_{12}^{(1)}\p_{r} V_2+ A_{21}^{(1)}\p_{\theta} V_1+e_1(r)V_1+ \tilde{e}_2(r) V_2\\
=\mathbb{F}(\bar{{\bf U}}^{(1)},\bar{{\bf U}}^{(2)}, \bar{B},\bar{A}),\\
\frac{1}{r}\p_{\theta}V_1-\frac1r\p_r(rV_2)=F_2(\bar{{\bf U}}^{(1)},\bar{B},\bar{A})-F_2(\bar{{\bf U}}^{(2)},\bar{B},\bar{A}),\\
V_1(r_0,\theta)-l_0 V_2(r_0,\theta)=0,\\
V_2(r_1,\theta)=0,
\end{cases}
\end{eqnarray}
where $A_{ij}^{(k)}=A_{ij}(\bar{{\bf U}}^k, \bar{B}, \bar{A})$ for any $i,j, k=1,2$ and
\be\no
&&\mathbb{F}(\bar{{\bf U}}^{(1)},\bar{{\bf U}}^{(2)}, \bar{B},\bar{A})=F_1(\bar{{\bf U}}^{(1)},\bar{B},\bar{A})-F_1(\bar{{\bf U}}^{(2)},\bar{B},\bar{A})-(A_{11}^{(1)}-A_{11}^{(2)})\p_r \hat{U}_2^{(2)}\\\no
&&\quad-r(A_{22}^{(1)}-A_{22}^{(2)})\p_{\theta} \hat{U}_2^{(2)}- r (A_{12}^{(1)}-A_{12}^{(2)})\p_r \hat{U}_2^{(2)}- (A_{21}^{(1)}-A_{21}^{(2)})\p_{\theta} \hat{U}_1^{(2)}
\ee

As above, decompose $V_1$ and $V_2$ as
\be\no
V_1=\frac{1}{r}\p_{\theta} \phi_3+ \p_r \phi_4,\ \ \ V_2= -\p_r\phi_3 + \frac{1}{r}\p_{\theta}\phi_4-\frac{r_1}{r} d_1, \ \ d_1= -\frac{1}{2\pi}\int_0^{2\pi} \p_r\phi_3(r_1,\theta) d\theta,
\ee
where $\phi_3$ and $\phi_4$ solve the following boundary value problems respectively:
\be\no\begin{cases}
(\p_r^2+\frac{1}{r}\p_r+\frac{1}{r^2}\p_{\theta}^2)\phi_3= F_2(\bar{{\bf U}}^{(1)},\bar{B},\bar{A})-F_2(\bar{{\bf U}}^{(2)},\bar{B},\bar{A}),\  \text{in }\Omega,\\
\phi_3(r_0,\theta)=\phi_1(r_1,\theta)=0,\ \ \theta\in \mathbb{T}_{2\pi},
\end{cases}\ee
and
\be\no\begin{cases}
A_{11}^{(1)}\p_r^2 \phi_4+r A_{22}^{(1)}\p_{\theta}^2 \phi_4+ (A_{12}^{(1)}+ A_{21}^{(1)})\p_{r\theta}^2 \phi_4+e_1(r)\p_r\phi_4+ E_2(\bar{{\bf U}}^{(1)}) \p_{\theta}\phi_4\\
=\mathbb{F}(\bar{{\bf U}}^{(1)},\bar{{\bf U}}^{(2)},\bar{B},\bar{A})- A_{11}^{(1)}\p_r(\frac{1}{r}\p_{\theta}\phi_3) +r A_{22}^{(1)}\p_{r\theta}^2\phi_3 + r A_{12}^{(1)}\p_{r}^2\phi_3 -\frac{A_{21}^{(1)}}{r} \p_{\theta}^2\phi_3\\\no
\quad- \frac{e_1(r)}{r} \p_{\theta}\phi_3 + (\tilde{e}_2(r)-A_{12}^{(1)})(\p_r\phi_3+\frac{r_1 d_1}{r}),\\
\p_r\phi_4(r_0,\theta)- l_0\frac{1}{r_0}\p_{\theta}\phi_4(r_0,\theta)=-\frac{l_0r_1 d_1}{r_0}-\frac{1}{r_0}\p_{\theta} \phi_3(r_0,\theta)- l_0\p_r\phi_3(r_0,\theta),\\
\frac{1}{r_1}\p_{\theta} \phi_4(r_1,\theta)= d_1 + \p_r \phi_3(r_1,\theta).
\end{cases}\ee
Then combining the $H^2(\Omega)$ estimate of $\phi_3$ and the $H^1(\Omega)$ estimate of $\phi_4$ leads to
\begin{align}\label{v1}
\|(V_1,V_2)\|_{L^2(\Omega)}&\leq C_1(\delta_0+\delta_1) \|(\bar{V}_1,\bar{V}_2)\|_{L^2(\Omega)}.
\end{align}
Choose $\delta_0=\sqrt{\eps+\delta_1}$ small enough such that
\begin{align}\no
\|(V_1,V_2)\|_{L^2(\Omega)}\leq \frac12 \|(\bar{V}_1,\bar{V}_2)\|_{L^2(\Omega)}
\end{align}
Then $\mathcal{T}^{(\bar{B},\bar{A})}$ is a contractive mapping in $L^2(\Omega)$-norm and there exists a unique fixed point $\bar{{\bf U}}\in \mathcal{X}_1$ to $\mathcal{T}^{(\bar{B},\bar{A})}$.


In summary, we have shown that, for any fixed $(\bar{B},\bar{A})\in\mathcal{X}_2$, the following problem has a unique solution $\bar{{\bf U}}\in\mathcal{X}_1$,
\begin{eqnarray}\label{374}
\begin{cases}
(c^2(\bar{B},|\bar{{\bf U}}|^2)-\bar{U}_1^2)\p_r \bar{U}_1+(c^2(\bar{B},|\bar{{\bf U}}|^2)-\bar{U}_2^2)\frac1r\p_\th \bar{U}_2\\-\bar{U}_1\bar{U}_2(\p_r\bar{U}_2+\frac1r\p_\th \bar{U}_1)+\frac1r c^2(\bar{B},|\bar{{\bf U}}|^2) \bar{U}_1\\
=-(\bar{U}_1\p_r+\bar{U}_2\frac1r\p_\th)\bar{B}+\frac{c^2(\bar{B},|\bar{{\bf U}}|^2)}{(\ga-1)\bar{A}}(\bar{U}_1\p_r+\bar{U}_2\frac1r\p_\th)\bar{A},\\
\frac{1}{r}\bar{U}_2(\p_\th \bar{U}_1-\p_r(r\bar{U}_2))=\frac{ \bar{B}-\frac{1}{2}|\bar{{\bf U}}|^2}{\bar{A}\gamma}\p_r \bar{A}-\p_r \bar{B},\\
\bar{U}_1(r_0,\th)-l_0 \bar{U}_2(r_0,\theta)=U_{b1}(r_0)-l_0U_{b2}(r_0)+\epsilon g_{0}(\theta)\in H^3(\mathbb{T}_{2\pi}),\\
\bar{U}_2(r_1,\th)=U_{20}+\epsilon g_{1}(\theta)\in C^{3,\alpha}(\mathbb{T}_{2\pi}).
\end{cases}
\end{eqnarray}

Note that when $\eps+\delta_1$ is suitably small, one may regard \eqref{374} as a uniformly first order elliptic system in $\tilde{\Omega}_{lue}$. Since $\bar{{\bf U}}\in \mathcal{X}_1$ and $\bar{B},\bar{A}\in\mathcal{X}_2$, so the coefficients in \eqref{374} belong to $H^3(\Omega)\subset C^{1,\alpha_1}(\overline{\Omega})$ for each $\alpha_1\in (0,1)$, the terms on the right hand side belong to $H^3(\Omega)$ and $\bar{U}_2(r_1,\th)\in C^{3,\alpha}(\mathbb{T}_{2\pi})$. Thus by standard interior and boundary regularity estimates to elliptic systems, one can improve the regularity of $\bar{{\bf U}}\in H^{4}(\Omega_{lue})\subset C^{2,\alpha}(\overline{\Omega_{lue}})$. This, together with the assumption $(\bar{B},\bar{A})\in C^{3,\alpha}(\overline{\Omega_{lue}})$, implies that the terms on the right hand side belong to $C^{2,\alpha}(\overline{\Omega_{lue}})$. The interior and boundary Schauder estimates to elliptic systems in $\Omega_{lue}$ yield that $\bar{{\bf U}}\in C^{3,\alpha}(\overline{\Omega_{ue}})$. In particular, $(\rho(\bar{{\bf U}},\bar{B},\bar{A})\bar{U}_1)(r_1,\cdot)\in C^{3,\alpha}(\mathbb{T}_{2\pi})$.


Next, for any $(\bar{B},\bar{A})\in\mathcal{X}_2$, we construct an operator $\mathcal{P}$: $(\bar{B},\bar{A})\in \mathcal{X}_2\mapsto(B,A)\in \mathcal{X}_2$, where $(B,A)=(B_0+\hat{B}, A_0+\hat{A})$ solves the following transport equations
\begin{eqnarray}\label{7step1}
\begin{cases}
\rho(\bar{\U},\bar{B},\bar{A})(r\bar{U}_1\partial_r +\bar{U}_2\partial_{\theta}) \hat{B}=0,\\
\rho(\bar{\U},\bar{B},\bar{A})(r\bar{U}_1\partial_r +\bar{U}_2\partial_{\theta})\hat{A}=0,\\
\hat{B}(r_1,\theta)=\epsilon B_1(\theta)\in C^{4,\alpha}(\mathbb{T}_{2\pi}),\\
\hat{A}(r_1,\theta)=\epsilon A_1(\theta)\in C^{4,\alpha}(\mathbb{T}_{2\pi}),
\end{cases}
\end{eqnarray}
here $\bar{\U}\in H^3(\Omega)\cap C^{2,\alpha}(\overline{\Omega_{lue}})\cap C^{3,\alpha}(\overline{\Omega_{ue}})$ is the unique fixed point of $\mathcal{T}^{(\bar{B},\bar{A})}$, that is, $\bar{{\bf U}}$ is the unique solution to \eqref{374}. By \eqref{374}, it holds
\be\no
\p_r(r\rho(\bar{\U},\bar{B},\bar{A})\bar{U}_1)+\p_\th(\rho(\bar{\U},\bar{B},\bar{A})\bar{U}_2)=0,
\ee
which enables one to define the stream function on $[r_0,r_1]\times \mathbb{R}$ as
\begin{align}\no
\psi(r,\th)=\int_0^{\th}r_1(\rho(\bar{\U},\bar{B},\bar{A})\bar{U}_1)(r_1,\tau)d\tau-\int_{r_1}^r(\rho(\bar{\U},\bar{B},\bar{A})\bar{U}_2)(\tau,\th)d\tau.
\end{align}
Note that the function $\psi$ defined above may not be periodic in $\theta$. However, it holds true that
\begin{align}\no\begin{cases}
\p_r\psi(r,\th)=-(\rho(\bar{\U},\bar{B},\bar{A})\bar{U}_2)(r,\th)\in H^3(\Omega)\cap C^{2,\alpha}(\overline{\Omega_{lue}})\cap C^{3,\alpha}(\overline{\Omega_{ue}}),\\
\p_\th\psi(r,\th)=r(\rho(\bar{\U},\bar{B},\bar{A})\bar{U}_1)(r,\th)\in H^3(\Omega)\cap C^{2,\alpha}(\overline{\Omega_{lue}})\cap C^{3,\alpha}(\overline{\Omega_{ue}}),\\
\psi(r,\th)\in H^4(\Omega)\cap C^{3,\alpha}(\overline{\Omega_{lue}})\cap C^{4,\alpha}(\overline{\Omega_{ue}}).
\end{cases}
\end{align}
Since the background transonic flow is symmetric in $r$ with $U_{b1}<0$, and $\bar{{\bf U}}\in \mathcal{X}_1$,  $\p_\th\psi(r_1,\th)=r_1(\rho(\bar{\U},\bar{B},\bar{A})\bar{U}_1)(r_1,\th)<0$, the inverse function of $\psi(r_1,\cdot)$: $\th\in\mathbb{R}\mapsto t\in\mathbb{R}$ is well-defined and is denoted by $\psi_{r_1}^{-1}$: $t\in\mathbb{R}\mapsto\th\in\mathbb{R}$.

Define the functions
\begin{align}\label{BA}
\hat{B}(r,\th)=\eps B_1(\psi_{r_1}^{-1}(\psi(r,\th))), \hat{A}(r,\th)=\eps A_1(\psi_{r_1}^{-1}(\psi(r,\th))).
\end{align}
We claim that $\hat{B}$ and $\hat{A}$ defined in \eqref{BA} are periodic in $\theta$ with period $2\pi$. Indeed, it follows from the definition that $\hat{B}(r,\th+2\pi)=\eps B_1(\psi_{r_1}^{-1}\circ\psi(r,\th+2\pi))$ and  $\hat{B}(r,\th)=\eps B_1(\psi_{r_1}^{-1}\circ\psi(r,\th))$. Denote $\beta_0=\psi_{r_1}^{-1}\circ\psi(r,\th)$ and $\beta_1=\psi_{r_1}^{-1}\circ\psi(r,\th+2\pi)$. It suffices to show that $\beta_0+2\pi=\beta_1$. Since $\rho(\bar{\U},\bar{B},\bar{A})(\bar{U}_1, \bar{U}_2)$ is periodic in $\theta$ with period $2\pi$,
\be\no
\psi_{r_1}(\beta_1)&=&\psi(r,\theta+2\pi)= \psi(r,\theta)+ \int_{\theta}^{\theta+2\pi} \rho(\bar{\U},\bar{B},\bar{A})\bar{U}_1(r_1,\tau) d\tau\\\no
&=& \psi_{r_1}(\beta_0)+ \int_{\beta_1-2\pi}^{\beta_1} \rho(\bar{\U},\bar{B},\bar{A})\bar{U}_1(r_1,\tau) d\tau
\ee
and noting that $\psi_{r_1}(\beta_1)=\int_{0}^{\beta_1} \rho(\bar{\U},\bar{B},\bar{A})\bar{U}_1(r_1,\tau) d\tau$, thus
\be\no
\psi_{r_1}(\beta_0)= \psi_{r_1}(\beta_1)- \int_{\beta_1-2\pi}^{\beta_1} \rho(\bar{\U},\bar{B},\bar{A})\bar{U}_1(r_1,\tau) d\tau=\psi_{r_1}(\beta_1-2\pi).
\ee
By monotonicity of $\psi_{r_1}(\cdot)$, $\beta_1-2\pi=\beta_0$. It is easy to verify that the functions defined in \eqref{BA} yield the unique solution to \eqref{7step1}. Since $(\rho(\bar{\U},\bar{B},\bar{A})\bar{U}_1)(r_1,\cdot)\in C^{3,\alpha}(\mathbb{T}_{2\pi})$, then $\psi_{r_1}^{-1}\in C^{4,\alpha}(\mathbb{R})$ and one has
\begin{align}\no
\|(\hat{B},\hat{A})\|_{H^4(\Omega)}+\|(\hat{B},\hat{A})\|_{C^{3,\alpha}(\overline{\Omega_{lue}})}+\|(\hat{B},\hat{A})\|_{C^{4,\alpha}(\overline{\Omega_{ue}})}\leq C_2\eps.
\end{align}
Therefore, $\mathcal{P}$ is well-defined if one selects $\delta_1=\sqrt{\eps}$ and $\sqrt{\eps}\leq \frac{1}{C_2}$.

It remains to show that the mapping $\mathcal{P}$ is contractive in a low order norm for suitably small $\eps$. Let $(B^{(i)},A^{(i)})=\mathcal{P}(\bar{B}^{(i)},\bar{A}^{(i)}) (i=1,2)$ for any $(\bar{B}^{(i)},\bar{A}^{(i)}) \in \mathcal{X}_2, (i=1,2)$ and denote $\bar{B}^{(1)}-\bar{B}^{(2)}$ by $\bar{B}_d$, $\bar{A}^{(1)}-\bar{A}^{(2)}$ by $\bar{A}_d$, $B^{(1)}-B^{(2)}$ by $B_d$ and $A^{(1)}-A^{(2)}$ by $A_d$, respectively. Then, it follows from \eqref{BA} that
\begin{align}\no
\begin{cases}
\hat{B}^{(i)}(r,\th)=\eps B_1\circ(\psi^{(i)}_{r_1})^{-1}\circ \psi^{(i)}(r,\th),\\
\hat{A}^{(i)}(r,\th)=\eps A_1\circ(\psi^{(i)}_{r_1})^{-1}\circ \psi^{(i)}(r,\th),
\end{cases}
\end{align}
where
\begin{align}\no
\psi^{(i)}(r,\th)=\int_0^{\th}r_1(\rho(\bar{\U}^{(i)},\bar{B}^{(i)},\bar{A}^{(i)})\bar{U}_1^{(i)})(r_1,\tau)d\tau-\int_{r_1}^r(\rho(\bar{\U}^{(i)},\bar{B}^{(i)},\bar{A}^{(i)})\bar{U}^{(i)}_2)(\tau,\th)d\tau,
\end{align}
$(\psi_{r_1}^{(i)})^{-1}$: $t\in\mathbb{R}\mapsto\th\in\mathbb{R}$ is the inverse function of $\psi^{(i)}(r_1,\cdot)$: $\th\in\mathbb{R}\mapsto t\in\mathbb{R}$ and
$\bar{\U}^{(i)}$ is the unique fixed point of $\mathcal{T}^{(\bar{B}^{(i)},\bar{A}^{(i)})}$, $i=1,2$.
Thus,
\begin{align}\no
|B_d|=|\hat{B}^{(1)}-\hat{B}^{(2)}|\leq\eps \|B'_1\|_{L^\infty(\mathbb{T}_{2\pi})}|\beta^{(1)}(r,\th)-\beta^{(2)}(r,\th)|
\end{align}
where $\beta^{(i)}(r,\th)=(\psi^{(i)}_{r_1})^{-1}\circ \psi^{(i)}(r,\th)\in[0,2\pi]$. It follows from the definitions that
\begin{align*}
&\int_{\beta^{(2)}(r,\th)}^{\beta^{(1)}(r,\th)}r_1(\rho(\bar{\U}^{(1)},\bar{B}^{(1)},\bar{A}^{(1)})\bar{U}_1^{(1)})(r_1,\tau)d\tau
=\psi^{(1)}(r,\th)-\psi^{(2)}(r,\th)\\
&\quad\quad\quad-\int_0^{\beta^{(2)}(r,\th)}r_1\{\rho(\bar{\U}^{(1)},\bar{B}^{(1)},\bar{A}^{(1)})\bar{U}_1^{(1)}-
\rho(\bar{\U}^{(2)},\bar{B}^{(2)},\bar{A}^{(2)})\bar{U}_1^{(2)}\}(r_1,\tau)d\tau,
\end{align*}
which implies
\begin{align*}
&\underline{m}^{(1)}|\beta^{(1)}(r,\th)-\beta^{(2)}(r,\th)|\\
\leq& |\psi^{(1)}(r,\th)-\psi^{(2)}(r,\th)|+r_1\int_0^{2\pi}|\rho(\bar{\U}^{(1)},\bar{B}^{(1)},\bar{A}^{(1)})\bar{U}_1^{(1)}-\rho(\bar{\U}^{(2)},\bar{B}^{(2)},\bar{A}^{(2)})\bar{U}_1^{(2)}|(r_1,\tau)d\tau
\end{align*}
with $\underline{m}^{(i)}:=r_1\min_{\th\in[0,2\pi]}(-\rho(\bar{\U}^{(i)},\bar{B}^{(i)},\bar{A}^{(i)})\bar{U}_1^{(i)})(r_1,\th)>0$. Noting that $\bar{B}_d(r_1,\th)=\bar{A}_d(r_1,\th)\equiv 0$, one has
\begin{align}\no
\|B_d\|_{L^2(\Omega)}\leq C_3\eps\bigg(\|(\bar{\U}^{(1)}-\bar{\U}^{(2)},\bar{B}_d,\bar{A}_d)\|_{L^2(\Omega)}+\|(\bar{\U}^{(1)}-\bar{\U}^{(2)})(r_1,\cdot)\|_{L^2(\mathbb{T}_{2\pi})}\bigg).
\end{align}
Since
\begin{align*}
|\p_r B_d|=&\eps|B_1'(\beta^{(1)}(r,\th))\p \beta^{(1)}(r,\th)-B_1'(\beta^{(2)}(r,\th))\p_r \beta^{(2)}(r,\th)|\\\no
=&\eps|\big(B_1'(\beta^{(1)}(r,\th))-B_1'(\beta^{(2)}(r,\th))\big)\p_r \beta^{(1)}+ B_1'(\beta^{(2)}(r,\th))(\p_r \beta^{(1)}-\p_r \beta^{(2)})|\\\no
\leq&\eps \|B''_1\|_{L^\infty(\mathbb{T}_{2\pi})}|\beta^{(1)}(r,\th)-\beta^{(2)}(r,\th)|\fr{1}{\underline{m}^{(1)}}\|\nabla \psi^{(1)}(r,\th)\|_{L^\infty(\Omega)}\\\no
&+\eps \|B'_1\|_{L^\infty(\mathbb{T}_{2\pi})}\fr{\|\nabla \psi^{(1)}(r,\th)\|_{L^\infty(\Omega)}}{\underline{m}^{(1)}\underline{m}^{(2)}}\bigg|(\rho(\bar{\U}^{(1)},\bar{B}^{(1)},\bar{A}^{(1)})\bar{U}_1^{(1)})(r_1,\psi^{(1)})\\
&\quad-(\rho(\bar{\U}^{(2)},\bar{B}^{(2)},\bar{A}^{(2)})\bar{U}_1^{(2)})(r_1,\psi^{(2)})\bigg|\\
&+\eps \|B'_1\|_{L^\infty(\mathbb{T}_{2\pi})}\fr{1}{\underline{m}^{(2)}} \bigg|\nabla \psi^{(1)}(r,\th)-\nabla \psi^{(2)}(r,\th)\bigg|,
\end{align*}
and similar computations are valid for $\p_\th B_d$. Then one has
\begin{align}\no
\|\nabla B_d\|_{L^2(\Omega)}\leq C_3\eps\bigg(\|(\bar{\U}^{(1)}-\bar{\U}^{(2)},\bar{B}_d,\bar{A}_d)\|_{L^2(\Omega)}+\|(\bar{\U}^{(1)}-\bar{\U}^{(2)})(r_1,\cdot)\|_{L^2(\mathbb{T}_{2\pi})}\bigg)
\end{align}
with $C_3$ independent of $\bar{\U}^{(i)},\bar{B}^{(i)},\bar{A}^{(i)}$, $i=1,2$.
Same estimate holds for $A_d$. Therefore,
\begin{align}\label{dd10}
\|(B_d,A_d)\|_{H^1(\Omega)}\leq C_3\eps\bigg(\|(\bar{\U}^{(1)}-\bar{\U}^{(2)},\bar{B}_d,\bar{A}_d)\|_{L^2(\Omega)}+\|(\bar{\U}^{(1)}-\bar{\U}^{(2)})(r_1,\cdot)\|_{L^2(\mathbb{T}_{2\pi})}\bigg).
\end{align}

We further claim that
\begin{eqnarray}\label{dd11}\|(\bar{\U}^{(1)}-\bar{\U}^{(2)})\|_{L^2(\Omega)}+\|(\bar{\U}^{(1)}-\bar{\U}^{(2)})(r_1,\cdot)\|_{L^2(\mathbb{T}_{2\pi})}\leq C_4\|(\bar{B}_d,\bar{A}_d)\|_{H^1(\Omega)}.
\end{eqnarray}

Indeed, set $\bar{U}_j^{(1)}-\bar{U}_j^{(2)}$ by $\bar{U}_{dj}$ for $j=1,2$ respectively. It follows from \eqref{374} that
\begin{eqnarray}\label{contract}
\begin{cases}
A_{11}(\bar{U}^{(1)},\bar{B}^{(1)})\p_r \bar{U}_{d1}+rA_{22}(\bar{U}^{(1)},\bar{B}^{(1)})\p_\th \bar{U}_{d2}
+rA_{12}(\bar{U}^{(1)},\bar{B}^{(1)})\p_r\bar{U}_{d2}\\ \quad\quad+A_{21}(\bar{U}^{(1)},\bar{B}^{(1)})\p_\th \bar{U}_{d1}+e_1(r)\bar{U}_{d1}+\tilde{e}_2(r)\bar{U}_{d2}=R_1,\\
\frac{1}{r}(\p_\th \bar{U}_{d1}-\p_r(r\bar{U}_{d2}))
=R_2,\\
\bar{U}_{d1}(r_0,\th)-l_0 \bar{U}_{d2}(r_0,\theta)=0,\\
\bar{U}_{d2}(r_1,\th)=0,
\end{cases}
\end{eqnarray}

where $R_1$ and $R_2$ are two quantities which satisfy
\begin{align}\no
\begin{cases}
\|R_1\|_{L^2(\Omega)}\leq C\bigg(\|(\bar{B}_d,\bar{A}_d)\|_{L^2(\Omega)}+\delta_0\|(\bar{\U}^{(1)}-\bar{\U}^{(2)})\|_{L^2(\Omega)}\bigg),\\
\|R_2\|_{L^2(\Omega)}\leq C\bigg(\|(\bar{B}_d,\bar{A}_d)\|_{H^1(\Omega)}+\delta_1\|(\bar{\U}^{(1)}-\bar{\U}^{(2)})\|_{L^2(\Omega)}\bigg).
\end{cases}
\end{align}


Similar arguments as for \eqref{v1} yield
\be\no
\|(\bar{\U}^{(1)}-\bar{\U}^{(2)})\|_{L^2(\Omega)}\leq C_4\bigg(\|(\bar{B}_d,\bar{A}_d)\|_{H^1(\Omega)}+(\delta_0+\delta_1)\|(\bar{\U}^{(1)}-\bar{\U}^{(2)})\|_{L^2(\Omega)}\bigg).
\ee
Since \eqref{contract} is uniformly elliptic in $\overline{\Omega_{lue}}$, the interior and boundary $H^1$ estimates for elliptic systems yield that
\be\no
\|(\bar{\U}^{(1)}-\bar{\U}^{(2)})\|_{H^1(\overline{\Omega_{ue}})}\leq C_4\bigg(\|(\bar{B}_d,\bar{A}_d)\|_{H^1(\Omega)}+(\delta_0+\delta_1)\|(\bar{\U}^{(1)}-\bar{\U}^{(2)})\|_{L^2(\Omega)}\bigg),
\ee
which further implies, by the trace Theorem, that
\be\no
\|(\bar{\U}^{(1)}-\bar{\U}^{(2)})(r_1,\cdot)\|_{L^2(\mathbb{T}_{2\pi})}\leq C_4\bigg(\|(\bar{B}_d,\bar{A}_d)\|_{H^1(\Omega)}+(\delta_0+\delta_1)\|(\bar{\U}^{(1)}-\bar{\U}^{(2)})\|_{L^2(\Omega)}\bigg).
\ee
Choosing $\delta_0+\delta_1=\sqrt{\eps+\delta_1}+\delta_1=\sqrt{\eps+\sqrt{\eps}}+\sqrt{\eps}$ small enough such that $C_4(\delta_0+\delta_1)<1/2$, one obtains \eqref{dd11}.

Combining \eqref{dd10} and \eqref{dd11}, we obtain finally that
\begin{align}\no
\|(B_d,A_d)\|_{H^1(\Omega)}&\leq C_5\eps \|(\bar{B}_d,\bar{A}_d)\|_{H^1(\Omega)}\leq \frac12 \|(\bar{B}_d,\bar{A}_d)\|_{H^1(\Omega)},
\end{align}
provided that $0<\epsilon \leq \frac{1}{2C_5}$. Hence $\mathcal{P}$ is a contractive mapping in $H^1(\Omega)$-norm and there exists a unique fixed point $(B,A)\in \mathcal{X}_2$. Denote the fixed point of the mapping $\mathcal{T}^{B,A}$ in $\mathcal{X}_1$ by ${\bf U}$, then $({\bf U}, B, A)$ is a solution to the 2-D steady Euler system \eqref{2d} with boundary conditions \eqref{2dboundary1}-\eqref{2dboundary4}, which also satisfies the estimate \eqref{2d10}. The uniqueness can be proved by a similar argument as for the contraction of the two mappings $\mathcal{T}^{B,A}$ and $\mathcal{P}$. The properties of the sonic surface can be proved as in Theorem \ref{main1}. The proof of Theorem \ref{2dmain} is completed.


\section{Smooth axi-symmetric transonic flows with small nonzero vorticity} \label{cylinder}\noindent

In this section, we prove Theorem \ref{cylind-theorem1}. As in the proof of Theorem \ref{2dmain}, we employ the deformation-curl decomposition in \cite{WengXin19,weng2019} and rewrite the system \eqref{ProblemImm} as
\begin{eqnarray}\label{rotational}
\begin{cases}
(c^2(B, |{\bf U}|^2)-U_1^2)\p_rU_1+(c^2(B, |{\bf U}|^2)-U_3^2)\p_3U_3-U_1U_3(\p_3U_1+\p_rU_3)+\frac{c^2(B, |{\bf U}|^2)+U_2^2}{r}U_1=0,\\
(U_1\partial_r +U_3 \partial_3) (rU_2)=0,\\
U_1(\p_rU_3-\p_3U_1)=-\p_3B+U_2\p_3U_2+\frac{1}{\ga-1}\rho^{\ga-1}\p_2A,\\
(U_1\partial_r+U_3 \partial_3) B=0,\\
(U_1\partial_r+U_3 \partial_3) A=0,
\end{cases}
\end{eqnarray}
with boundary conditions on $r=r_1$:
\begin{align}
U_3=\eps q_3(x_3),\ \  U_2= U_{b2}+\eps q_2(x_3),\ \ B=B_0+\eps \tilde{B}_1(x_3),\ \ A=A_0+\eps\tilde{A}_1(x_3),
\end{align}
and on $r=r_0$:
\begin{align}
U_1=U_{b1}+\eps q_1(x_3).
\end{align}
Here, $q_i$, $i=1,2,3$, $\tilde{B}_1$, $\tilde{A}_1\in C_c^{2,\alpha}(\mathbb{R})$.

As before, set
\begin{eqnarray}\no
\hat{U}_1= U_1- U_{b1},\ \ \hat{U}_2= U_2- U_{b2},\ \ \hat{U}_3= U_3,\ \ \hat{B}= B-B_0,\ \ \hat{A}= A-A_0.
\end{eqnarray}

Then ${\bf \hat{U}}$, $\hat{B}$ and $\hat{A}$ satisfy
\begin{eqnarray}\label{cylind1}
\begin{cases}
A_{b11}\p_r\hat{U}_1+A_{b33}\p_3\hat{U}_3+e_1(r)\hat{U}_1=G_1({\bf U},B,A),\\
\p_r\hat{U}_3-\p_3\hat{U}_1=G_2({\bf U},B,A),\\
(U_1\partial_r +U_3 \partial_3) (r\hat{U}_2)=0,\\
(U_1\partial_r+U_3 \partial_3) \hat{B}=0,\\
(U_1\partial_r+U_3 \partial_3) \hat{A}=0,
\end{cases}
\end{eqnarray}
where $e_1(r)$ is given as Section \ref{key}, $A_{b11}(r)=c^2(\rho_b)-U_{b1}^2,\ A_{b33}(r)=c^2(\rho_b)$ and
\begin{eqnarray}\label{cylind2}
\begin{cases}
G_1({\bf U},B,A)=-(\ga-1)(U'_{b1}+\frac{1}{r}U_{b1})(\hat{B}-\frac12\hat{U}_3^2)+\bigg((\ga-1)U'_{b1}U_{b2}+\frac{\ga-3}{r}U_{b1}U_{b2}\bigg)\hat{U}_2\\
+(\frac{\ga+1}{2}U'_{b1}+\frac{\ga-1}{2r}U_{b1})\hat{U}_1^2+(\frac{\ga-1}{2}U'_{b1}+\frac{\ga-3}{2r}U_{b1})\hat{U}_2^2-\frac{1}{r}(c^2+U_2^2-c_b^2-U_{b2}^2)\hat{U}_1\\
-(c^2-U_1^2-c_b^2+U_{b1}^2)\p_r\hat{U}_1-(c^2-c_b^2-\hat{U}_3^2)\p_3\hat{U}_3+U_1\hat{U}_3(\p_r\hat{U}_3+\p_3\hat{U}_1),\\
G_2({\bf U},B,A)=\dfrac{1}{U_1}\bigg(-\p_3\hat{B}+U_2\p_3\hat{U}_2+\frac{1}{\ga-1}\rho^{\ga-1}\p_3\hat{A}\bigg).
\end{cases}
\end{eqnarray}
In contracts to the mixed type system \eqref{aa}, the left hand side of the system \eqref{cylind1} is a decoupled hyperbolic-elliptic system, $rU_2$, $B$ and $A$ satisfy hyperbolic equations, while $U_1$ and $U_3$ satisfy a first order elliptic system. By resolving the hyperbolic quantities first, it is easy to see that $\hat{U}_2, \hat{B}$ and $\hat{A}$ are of same order as $O(\epsilon)$. The terms involving derivatives $\nabla \hat{U}_j$ for $j=1,3$ in $G_1$ contain also a small factor, thus the left hand side of the first two equations in \eqref{cylind1} are the principal part.

Now we start to prove Theorem \ref{cylind-theorem1}. Define the solution space as
\begin{eqnarray}\nonumber
\mathcal{X}=\left\{({\bf U},B,A)(r,x_3) \in C^{2,\alpha}(\overline{\mathbb{D}}): \|(\hat{{\bf U}},\hat{B},\hat{A})\|_{C^{2,\alpha}(\overline{\mathbb{D}})}\leq \delta_0,
\displaystyle\lim_{x_3\to \pm\infty}(\hat{{\bf U}},\hat{B},\hat{A})(r,x_3)=0\right\},
\end{eqnarray}
with $\delta_0>0$ to be specified later. For any $(\bar{{\bf U}},\bar{B},\bar{A})\in \mathcal{X}$, we will construct an operator $\mathcal{T}$: $(\bar{{\bf U}},\bar{B},\bar{A})\in \mathcal{X}\mapsto({\bf U},B,A)\in \mathcal{X}$, with $({\bf U},B,A)$ to be obtained by the following steps.

First one obtains $(U_2,B,A)$ by solving the following hyperbolic problems:
\begin{eqnarray}\label{step1}
\begin{cases}
(\bar{U}_1\partial_r +\bar{U}_3 \partial_3) (r\hat{U}_2, \hat{B},\hat{A})=0,\\
(\hat{U}_2,\hat{B},\hat{A})(r_1,x_3)=(\epsilon q_{2}(x_3),\epsilon \tilde{B}_1(x_3),\epsilon \tilde{A}_1(x_3)).
\end{cases}
\end{eqnarray}
Since $\bar{{\bf U}}\in \mathcal{X}$ and $\bar{U}_1>0$ and $\|\bar{U}_3\|_{C^{2,\alpha}(\overline{\mathbb{D}})}\leq \delta_0$, the above transport equations can be solved by the characteristics method with the following estimates
\begin{align}\label{cylind3}
\|(\hat{U}_2,\hat{B},\hat{A})\|_{C^{2,\alpha}(\overline{\mathbb{D}})}\leq C_6\epsilon,
\end{align}
where $C_6$ depends only on the background solution and the boundary datum. Moreover, since $(q_2,\tilde{B}_1, \tilde{A}_1)$ has compact support, it is easy to see that $(\hat{U}_2, \hat{B},\hat{A})$ also has compact support and
\be\label{cylin31}
\displaystyle\lim_{x_3\to \pm \infty} (\hat{U}_2, \hat{B},\hat{A})(r,x_3)=0, \ \lim_{x_3\to \pm \infty} \nabla_{r,x_3}(\hat{U}_2, \hat{B},\hat{A})(r,x_3)=0.
\ee

Next we will solve the following boundary value problem for a linear first order elliptic system to obtain $(U_1,U_3)$.
\begin{eqnarray}\label{step2}
\begin{cases}
A_{b11}(r)\p_r\hat{U}_1+A_{b33}(r)\p_3\hat{U}_3+e_1(r)\hat{U}_1=G_1(\bar{U}_1,U_2,\bar{U}_3,B,A),\\
\p_r\hat{U}_3-\p_3\hat{U}_1=G_2(\bar{U}_1,U_2,\bar{U}_3,B,A),\\
\hat{U}_1(r_0,x_3)=\epsilon q_{1}(x_3),\\
\hat{U}_3(r_1,x_3)=\epsilon q_{3}(x_3),\\
\displaystyle\lim_{x_3\to\pm \infty} \hat{U}_1(r,x_3)=0.
\end{cases}
\end{eqnarray}
With the formula of $G_1$ and $G_2$ in \eqref{cylind2} and the estimates \eqref{cylind3}, it can be verified directly that
\begin{align}\label{cylind4}
\begin{cases}
\|G_1(\bar{U}_1,U_2,\bar{U}_3,B,A)\|_{C^{1,\alpha}(\overline{\mathbb{D}})}\leq C_6(\epsilon+\delta_0^2),\\
\|G_2(\bar{U}_1,U_2,\bar{U}_3,B,A)\|_{C^{1,\alpha}(\overline{\mathbb{D}})}\leq C_6\epsilon.,\\
\displaystyle\lim_{x_3\to\pm\infty} G_1(\bar{U}_1,U_2,\bar{U}_3,B,A)=0,\ \lim_{x_3\to\pm\infty}  G_2(\bar{U}_1,U_2,\bar{U}_3,B,A)=0.
\end{cases}
\end{align}
It is easy to show that there exists a unique solution $\phi_1(r,x_3)$ to the following problem
\begin{eqnarray}\label{cylind5}
\begin{cases}
(\p_r^2+\p_3^2)\phi_1=G_2(\bar{U}_1,U_2,\bar{U}_3,B,A)\in C^{1,\alpha}(\overline{\mathbb{D}}),\\
\phi_1(r_0,x_3)=\p_r\phi_1(r_1,x_3)=0,\\
\displaystyle\lim_{x_3\to\pm\infty} \p_3 \phi_1(r,x_3)= 0.
\end{cases}
\end{eqnarray}
Moreover, $\phi_1(r,x_3)\in C^{3,\alpha}(\overline{\mathbb{D}})$ with the property that
\begin{eqnarray}\label{cylind6}
&&\|\phi_1\|_{C^{3,\alpha}(\overline{\mathbb{D}})}\leq C\|G_2(\bar{U}_1,U_2,\bar{U}_3,B,A)\|_{C^{1,\alpha}(\overline{\mathbb{D}})}\leq C_6\epsilon,\\\label{cylind7}
&&\displaystyle\lim_{x_3\to\pm\infty} (\nabla_{r,x_3}\phi_1, \nabla_{r,x_3}^2\phi_1)(r,x_3)=0.
\end{eqnarray}
Define $V_1= \hat{U}_1+ \p_3 \phi_1$ and $V_3= \hat{U}_3- \partial_r \phi_1$. Then
\begin{eqnarray}\label{step21}
\begin{cases}
A_{b11}(r)\p_r V_1+A_{b33}(r)\p_3V_3+e_1(r)V_1=G_3(\bar{U}_1,U_2,\bar{U}_3,B,A),\\
\p_r V_3-\p_3 V_1=0,\\
V_1(r_0,x_3)=\epsilon q_{1}(x_3),\\
V_3(r_1,x_3)=\epsilon q_{3}(x_3),\\
\displaystyle\lim_{x_3\to \pm\infty} V_1(r,x_3)=0,
\end{cases}\end{eqnarray}
where
\be\no
G_3(\bar{U}_1,U_2,\bar{U}_3,B,A)= G_1(\bar{U}_1,U_2,\bar{U}_3,B,A)- U_{b1}^2(r)\p_{rx_3}^2 \phi_1 + e_1(r)\p_3 \phi_1.
\ee

The second equation in \eqref{step21} implies that there exists a potential function $\phi(r,x_3)$ such that $V_1=\p_r \phi, V_3=\partial_3 \phi$ and $\phi$ should satisfy the following second-order elliptic equation in $\mathbb{D}$:
\begin{eqnarray}\label{step23}
\begin{cases}
A_{b11}(r)\p_r^2\phi+A_{b33}(r)\p_3^2\phi+e_1(r)\p_r\phi=G_3(\bar{U}_1,U_2,\bar{U}_3,B,A),\  &\text{in }\mathbb{D},\\
\p_r\phi(r_0,x_3)=\epsilon q_1(x_3),\ &\forall x_3\in\mathbb{R},\\
\p_3\phi(r_1,x_3)=\eps q_3(x_3),\ &\forall x_3\in\mathbb{R},\\
\displaystyle\lim_{x_3\to \pm\infty}\p_r\phi(r,x_3)=0,\ \ \phi(r_1,0)=0.
\end{cases}
\end{eqnarray}

To prove the existence and uniqueness of smooth solution to \eqref{step23}, one may first consider the problem in a truncated domain
\begin{eqnarray}\label{step24}
\begin{cases}
A_{b11}(r)\p_r^2\phi_n+A_{b33}(r)\p_3^2\phi_n+e_1(r)\p_r\phi_n=G_3,\ &\text{in }\mathbb{D}_n:=(r_0,r_1)\times (-n,n),\\
\p_r\phi_n(r_0,x_3)=\epsilon q_1(x_3),\ &\forall x_3\in [-n,n]\\
\phi_n(r_1,x_3)=\epsilon \int_0^{x_3} q_3(s) ds,\ &\forall x_3\in [-n,n],\\
\phi_n(r, n)=\epsilon \int_0^{n} q_3(s) ds,\ &\forall r\in (r_0,r_1),\\
\phi_n(r, -n)=\epsilon \int_0^{-n} q_3(s) ds,\ &\forall r\in (r_0,r_1).
\end{cases}
\end{eqnarray}

By Theorem 1 in \cite{lieberman86} and the remark after that theorem, there exists a unique solution $\phi_n\in C^2(\mathbb{D}_n)\cap C(\overline{\mathbb{D}_n})$ to \eqref{step24}. It remains to derive some uniform estimates in $n$ and take a limit to obtain a solution to \eqref{step23}. Define a barrier function $v(r,x_3)= m_1(\|G_3\|_{L^{\infty}}+\epsilon\|q_1\|_{L^{\infty}}) r$, with constant $m_1$ to be specified. Then $\phi_n-v$ satisfies
\begin{eqnarray}\label{step25}
\begin{cases}
A_{b11}(r)\p_r^2(\phi_n-v)+A_{b33}(r)\p_3^2(\phi_n-v)+e_1(r)\p_r(\phi_n-v)\\
\quad\quad\quad\quad=G_3-m_1 e_1(\|G_3\|_{L^{\infty}}+\epsilon\|q_1\|_{L^{\infty}}),\\
\p_r(\phi_n-v)(r_0,x_3)=\eps q_1(x_3)-m_1(\|G_3\|_{L^{\infty}}+\epsilon\|q_1\|_{L^{\infty}}),\\
(\phi_n-v)(r_1,x_3)=\epsilon \int_0^{x_3} q_3(s) ds-m_1 r_1(\|G_3\|_{L^{\infty}}+\epsilon\|q_1\|_{L^{\infty}}),\\
(\phi_n-v)(r,n)=\epsilon \int_0^{n} q_3(s) ds-m_1(\|G_3\|_{L^{\infty}}+\epsilon\|q_1\|_{L^{\infty}})r, \\
(\phi_n-v)(r,-n)=\epsilon \int_0^{-n} q_3(s) ds-m_1(\|G_3\|_{L^{\infty}}+\epsilon\|q_1\|_{L^{\infty}}) r.
\end{cases}
\end{eqnarray}
Since $e_1(r)>0$ for all $r\in [r_0,r_1]$, one may choose $m_1<0$ independent of $n$ such that
\begin{eqnarray}\nonumber
&&G_3-m_1 e_1(\|G_3\|_{L^{\infty}}+\epsilon\|q_1\|_{L^{\infty}})>0,\ \ \forall (r,x_3)\in \Omega_n; \\\nonumber
&&\epsilon q_1(x_3)-m_1(\|G_3\|_{L^{\infty}}+\epsilon\|q_1\|_{L^{\infty}})>0,\ \ \forall x_3\in [-n,n].
\end{eqnarray}
Therefore, the maximum principle shows that $\phi_n-v$ attains its maximum only on its boundary except $\{(r_0,x_3):x_3\in[-n,n]\}$ and thus
\begin{eqnarray}\nonumber
\phi_n\leq C_6(\|G_3\|_{L^{\infty}}+\epsilon\|q_1\|_{L^{\infty}}+ \epsilon \|q_3\|_{L^{1}(\mathbb{R})}),
\end{eqnarray}
where $C_6$ is a positive constant independent of $n$. Similarly, we can derive a lower bound for $\phi_n$. Thus the following uniform $L^{\infty}$ estimate holds:
\begin{eqnarray}\label{step25}
\|\phi_n\|_{L^{\infty}(\mathbb{D}_n)}\leq C_6(\|G_3\|_{L^{\infty}(\mathbb{D}_n)}+\displaystyle \epsilon\sum_{j=1,3}\|q_j\|_{L^{\infty}}+ \epsilon \|q_3\|_{L^{1}(\mathbb{R})}).
\end{eqnarray}

Utilizing Theorem 6.6 and Theorem 6.30 in \cite{gt}, for any compact domain $K\Subset [r_0,r_1]\times \mathbb{R}$, we have for any large $n$
\begin{eqnarray}\label{step26}
\|\phi_n\|_{C^{2,\alpha}(K)}\leq C_6(\|G_3\|_{C^{\alpha}(\overline{\mathbb{D}})}+\displaystyle\epsilon\sum_{j=1,3}\|q_j\|_{C^{1,\alpha}(\mathbb{R})}+ \epsilon \|q_3\|_{L^{1}(\mathbb{R})}),
\end{eqnarray}
where $C$ is independent of $n$. By a diagonal argument, one can extract a subsequence $\{\phi_{n_k}\}_{k=1}^{\infty}$ such that
\begin{eqnarray}\label{step27}
\phi_{n_k}\to \phi\ \ \text{in $C^{2,\beta}(K)$ for any compact subregion $K\Subset \overline{\mathbb{D}}$ and any $0<\beta<\alpha$}.
\end{eqnarray}
Hence $\phi$ admits the following estimate
\begin{eqnarray}\label{step28}
\|\phi\|_{C^{2,\alpha}(\overline{\mathbb{D}})}\leq C_6(\|G_3\|_{C^{\alpha}(\overline{\mathbb{D}})}+\displaystyle\sum_{j=1,3}\|q_j\|_{C^{1,\alpha}(\mathbb{R})}+ \epsilon \|q_3\|_{L^{1}(\mathbb{R})}).
\end{eqnarray}
Moreover, $\phi$ solves the following problem
\begin{eqnarray}\label{step29}
\begin{cases}
A_{b11}(r)\p_r^2\phi+A_{b33}(r)\p_3^2\phi+e_1(r)\p_r\phi=G_3(\bar{U}_1,U_2,\bar{U}_3,B,A),\  &\text{in }\mathbb{D},\\
\p_r\phi(r_0,x_3)=\epsilon q_1(x_3),\ &\forall x_3\in\mathbb{R},\\
\phi(r_1,x_3)=\epsilon \int_0^{x_3} q_3(s) ds,\ &\forall x_3\in\mathbb{R}.
\end{cases}
\end{eqnarray}

The asymptotic behavior of $\phi$ can be derived by a blowup argument. Define the functions $\psi_n(r,x_3)= \phi(r,x_3-n)$ for any $(r,x_3)\in [r_0,r_1]\times \mathbb{R}$. Then for any compact domain $K\subset [r_0,r_1]\times \mathbb{R}$, it follows from \eqref{step28} that
\begin{eqnarray}\nonumber
\|\psi_n\|_{C^{2,\alpha}(K)}\leq M_1:=C_6(\|G_3\|_{C^{\alpha}(\overline{\mathbb{D}})}+\displaystyle\sum_{j=1,3}\|q_j\|_{C^{1,\alpha}(\mathbb{R})}+ \epsilon \|q_3\|_{L^{1}(\mathbb{R})}).
\end{eqnarray}
Then it follows from Arzela-Ascoli theorem and a diagonal argument that there exists a subsequence $\psi_{n_k}$ such that
\begin{eqnarray}\nonumber
\psi_{n_k} \rightrightarrows \psi \ \text{as $n_k\to \infty$ in }\ C^{2,\beta}(K),
\end{eqnarray}
for any compact domain $K$ and any $\beta\in (0,\alpha)$. Therefore $\psi$ is bounded by $M_1$ in $[r_0,r_1]\times \mathbb{R}$. Note that $G_3(\bar{U}_1,U_2,\bar{U}_3,B,A)\to 0$ as $x_3\to +\infty$. Then $\psi$ should solve
\begin{eqnarray}\label{step29}
\begin{cases}
A_{b11}(r)\p_r^2\psi+A_{b33}(r)\p_3^2\psi+e_1(r)\p_r\psi=0,\  &(r,x_3)\in \mathbb{D},\\
\p_r\psi(r_0,x_3)=0,\ &\forall x_3\in\mathbb{R},\\
\psi(r_1,x_3)=\epsilon \int_0^{\infty} q_3(s) ds,\ &\forall x_3\in\mathbb{R}.
\end{cases}
\end{eqnarray}
We now prove $\psi\equiv \epsilon \int_0^{\infty} q_3(s) ds$. Set $\hat{\psi}=\psi-\epsilon \int_0^{\infty} q_3(s) ds$. For each $\eta>0$, choose a barrier function as $b(r,x_3)= \eta x_3^2- \mu \eta(r-r_1)$, where $\mu$ is any fixed constant larger than $\displaystyle\sup_{r\in [r_0,r_1]}\frac{2 A_{33}(r)}{e_1(r)}$. Then for large enough $m$, it holds that
\begin{eqnarray}\label{step30}
\begin{cases}
A_{b11}(r)\p_r^2(\hat{\psi}-b)+A_{b33}(r)\p_3^2(\hat{\psi}-b)+e_1(r)\p_r(\hat{\psi}-b)\\
\quad\quad\quad =\eta(\mu e_1(r)-2A_{b33}(r))\geq 0,\  &(r,x_3)\in (r_0,r_1)\times (-m,m),\\
\p_r(\hat{\psi}-b)(r_0,x_3)=\mu \eta> 0,\ &\forall x_3\in [-m,m],\\
(\hat{\psi}-b)(r_1,x_3)=-\eta x_3^2\leq 0,\ &\forall x_3\in [-m,m],\\
(\hat{\psi}-b)(r, \pm m)\leq M_1- \eta m^2+ \mu\eta (r-r_1)\leq 0, \ &\forall r\in (r_0,r_1).
\end{cases}
\end{eqnarray}
It follows from maximum principle that
\begin{eqnarray}\nonumber
-b(r,x_3)\leq \hat{\psi}(r,x_3)\leq b(r,x_3), \ \ \forall (r,x_3)\in [r_0,r_1]\times [-m,m].
\end{eqnarray}
For any fixed point $(r,x_3)$, letting $\eta\to 0$ shows that $\hat{\psi}\equiv 0$. Thus $\psi(r,x_3)\equiv \eps \int_0^{\infty} q_3(s) ds$, which implies that
\begin{eqnarray}\nonumber
\nabla \psi_{n_k}\rightrightarrows 0 \ \ \text{as $n_k\to \infty$ in }C^{1,\beta}(K).
\end{eqnarray}
Therefore $\nabla \phi(r,x_3)\to 0$ as $x_3\to +\infty$. Similarly, one can show that $\nabla \phi(r,x_3)\to 0$ as $x_3\to -\infty$. The existence and uniqueness of $\phi_1$ to \eqref{cylind5} can be proved in a similar way. Hence, one has shown that
\begin{eqnarray}\label{cylind11}
&&\|(\hat{U}_1,\hat{U}_3)\|_{C^{2,\alpha}(\overline{\mathbb{D}})}=\|(\p_r\phi-\p_3\phi_1,\p_3\phi+\p_r\phi_1)\|_{C^{2,\alpha}(\overline{\mathbb{D}})}\leq C_6(\epsilon+\delta_0^2),\\\label{cylind12}
&&\displaystyle \lim_{x_3\to \pm\infty} (\hat{U}_1,\hat{U}_3)(r,x_3)=0.
\end{eqnarray}
Set $\delta_0= 2C_* \eps$ and select a $0<\epsilon_0<\frac{1}{4C_1^2}$ small enough such that $C_6(\epsilon+\delta_0^2)\leq \delta_0$. Thus one has obtained a mapping $\mathcal{T}$ from $\mathcal{X}$ to itself. By a similar argument, one can prove that $\mathcal{T}$ is a contraction mapping in a weak norm. Thus there exists a unique fixed point $({\bf U}, B, A)\in \mathcal{X}$ to $\mathcal{T}$, which is the desired solution. The information about the sonic surface to the solution $({\bf U}, B, A)\in \mathcal{X}$ can be obtained in the same way as in Theorem \ref{main1}. The proof is completed.

{\bf Acknowledgement.}  Part of this work was done when Weng visited The Institute of Mathematical Sciences of The Chinese University of Hong Kong. He is grateful to the institute for providing nice research environment. Weng is partially supported by National Natural Science Foundation of China 11701431, 11971307, 12071359, the grant of Project of Thousand Youth Talents (No. 212100004). Xin is supported in part by the Zheng Ge Ru  Foundation, Hong Kong RGC Earmarked Research Grants CUHK-14305315, CUHK-14300917, CUHK-14302819 and CUHK-14302917, and by Guangdong Basic and Applied Basic Research Fundation 2020B1515310002.


\begin{thebibliography}{10}

\bibitem{bdxx19}
M. Bae, B. Duan, J. Xiao and C. Xie. {\it Structural stability of supersonic solutions to the Euler-Poisson system.} arXiv: 1602.01892v2. 10. Jan 2019.



\bibitem{bers54}
Lipman Bers, {\it Existence and uniqueness of a subsonic flow past a given profile.} Comm. Pure Appl. Math. 7 (1954), 441-504.

\bibitem{Bers1958}
Lipman Bers, {\it Mathematical Aspects of Subsonic and Transonic Gas
  Dynamics.}, Wiley, New York, 1958.

\bibitem{ccs06}
G. Chen, J. Chen and K. Song, {\it Transonic nozzle flows and free boundary problems for the full Euler equations,} J. Differential Equations, 229 (2006), no. 1, 92-120.

\bibitem{ccf07}
G. Chen, J. Chen and M. Feldman. {\it Transonic shocks and free boundary problems for the full Euler equations in infinite nozzles.} J. Math. Pures Appl. (9) 88 (2007), no. 2, 191-218.

\bibitem{cdsw07}
G. Chen, C. Dafermos, M. Slemrod, D. Wang, {\it On two-dimensional sonic-subsonic flow.} Comm. Math. Phys. 271 (2007), no. 3, 635-647.

\bibitem{chw16}
G. Chen, F. Huang, T. Wang, {\it Sonic-subsonic limit of approximate solutions to multidimensional steady Euler equations,} Arch. Rational Mech. Anal.,219 (2016), 719-740.

\bibitem{chen08}
S. Chen, {\it Transonic shocks in 3-D compressible flow passing a duct with a general section for Euler systems,} Trans. Amer. Math. Soc, 360 (2008), 5265-5289.

\bibitem{Courant1948}
R. Courant and K. O. Friedrichs, {\it Supersonic Flow and Shock Waves},
  Interscience Publishers, Inc., New York, 1948.

\bibitem{Cui2011}
D. Cui and J. Li, {\it On the existence and stability of 2-D perturbed
steady subsonic circulatory flows}, Science China Mathematics, 2011 Vol. 54, No. 7: 1421-1436.

\bibitem{fg57}
R. Finn and D. Gilbarg. {\it Asymptotic behavior and uniqueness of plane subsonic flows.} Comm. Pure Appl. Math., 10 (1957), 23-63.

\bibitem{Friedrichs1958}
K. O. Friedrichs, {\it Symmetric Positive Linear Differential Equations}, Comm. Pure Appl. Math. \textbf{XI} (1958), 333-418.

\bibitem{gt}
D. Gilbarg and N. Trudinger. {\it Elliptic Partial Differential Equations of Second Order.} 2nd Ed. Springer-Verlag: Berlin.

\bibitem{Gu1981}
C. Gu, {\it On partial differential equations of mixed type in n independent variables},
Comm. Pure Appl. Math. 34 (1981), 333-345.

\bibitem{hww11}
F. Huang, T. Wang, and Y. Wang, {\it On multi-dimensional sonic-subsonic
flow}, Acta Mathematica Scientia, 31 (2011), 2131-2140.

\bibitem{Kuzmin2002}
Alexander G. Kuz'min, {\it Boundary-Value Problems for Transonic Flow}, John
  Wiley ${\&}$ Sons, Ltd, West Sussex, 2002.

\bibitem{LP1969}
P. D. Lax, and R. W. Phillips, {\it Local boundary conditions for symmetric linear differential
operators}, Comm. Pure Appl. Math. 13 (1969), 427-455.

\bibitem{lieberman86}
G. M. Lieberman. {\it Mixed boundary value problems for elliptic and
parabolic differential equations of second order}, J. Math. Anal. Appl. 113, 1986,
pp. 422-440.

\bibitem{lxy09b}
J. Li, Z. Xin and H. Yin, {\it A free boundary value problem for the full Euler system and 2-D transonic shock in a large variable nozzle,} Math. Res. Lett., {\bf 16} (2009), 777-796.

\bibitem{lxy13}
J. Li, Z. Xin and H. Yin, {\it Transonic shocks for the full compressible Euler system in a general two-dimensional de Laval nozzle,} Arch. Ration. Mech. Anal. 207 (2013), no. 2, 533-581.

\bibitem{LMP07}
D. Lupo, C. Morawetz and K. Payne, {\it On Closed Boundary Value Problems for Equations
of Mixed Elliptic-Hyperbolic Type}, Comm. Pure Appl. Math., Vol. LX, 1319-1348 (2007).
%



\bibitem{Morawetz1956}
C. Morawetz. {\it On the non-existence of continuous transonic flows
past profiles I}, Comm. Pure Appl. Math. \textbf{9}
  (1956), no.~1, 45--68.

\bibitem{Morawetz1964}
C. Morawetz. {\it Non-existence of transonic flow past a profile},
  Comm. Pure Appl. Math. \textbf{17} (1964), 357--367.

\bibitem{Morawetz1958}
C. Morawetz. {\it A weak solution for a system of equations of elliptic-hyperbolic type},
Comm. Pure Appl. Math. 11 (1958), 315-331.

\bibitem{Morawetz1970}
C. Morawetz. {\it The Dirichlet problem for the tricomi equation}, Comm. Pure Appl.
Math. 23 (1970), 587-601.

\bibitem{Morawetz2004}
C. Morawetz. {\it Mixed equations and transonic flow.} J. Hyperbolic Differ. Equ. 1 (2004), no. 1,
1-26.

\bibitem{OD1955}
S. Ou and S. Ding. {\it Sur L'unicite du probleme de Tricomi de Leqaution de Chaplygin}, (Chinese, French Summary), Acta Math. Sin. Vol. 5, 1955, 697-703.

\bibitem{Protter1954}

M. H. Protter. {\it An existence theorem for the generalized Tricomi problem.} Duke Math. J. Vol. 21, 1954, 1-7.


\bibitem{Smirnov1978}
M. M. Smirnov, {\it Equations of Mixed Type}, Translations of Mathematical Monographs,
Vol. 51, Amer. Math. Soc., Providence, RI., 1978.

\bibitem{shiffman52}
M. Shiffman. {\it On the existence of subsonic flows of a compressible fluid.} J. Rational Mech. Anal., 1 (1952), 605-652.

\bibitem{WX2013}
C. Wang and Z. Xin. {\it On a degenerate free boundary problem and continuous subsonic-sonic flows in a convergent nozzle.} Arch. Ration. Mech. Anal. 208 (2013), no.3, 911-975.


\bibitem{WX2016}
C. Wang and Z. Xin. {\it On sonic curves of smooth subsonic-sonic and transonic flows.} SIAM J. Math. Anal. 48 (2016), no. 4. 2414-2453.

\bibitem{WX2019}
C. Wang and Z. Xin. {\it Smooth Transonic Flows of Meyer Type in
De Laval Nozzles}, Arch. Ration. Mech. Anal. \textbf{232}
  (2019), no.~3, 1597--1647.

\bibitem{WX2020}
C. Wang and Z. Xin. {\it Lipschitz Continuous Subsonic-sonic Flows in General Nozzles}, preprint 2019.

\bibitem{WengXin19}
S. Weng and Z. Xin. {\it A deformation-curl decomposition for three dimensional steady Euler equations (in Chinese).}
  Sci Sin Math, 2019, 49: 307-320, doi: 10.1360/N012018-00125.

\bibitem{weng2019}
S. Weng. {\it A deformation-curl-Poisson decomposition to the three dimensional steady Euler-Poisson system with applications.} J. Differential Equations 267 (2019), no. 11, 6574-6603.

\bibitem{WXY20a}
S. Weng, Z. Xin and H. Yuan. {\it Steady Compressible Radially Symmetric Flows in an Annulus}, preprint 2020.


\bibitem{xx07}
C. Xie and Z. Xin. {\it Global subsonic and subsonic-sonic flows through
infinitely long nozzles}, Indiana University Mathematical Journal, 56 (2007), 2991-3023.

\bibitem{xx10}
C. Xie and Z. Xin. {\it Global subsonic and subsonic-sonic flows through infinitely long axially symmetric nozzles}, Journal of Differential Equations \textbf{248} (2010), no.~11, 2657-2683.

\bibitem{xy05}
Z. Xin and H. Yin, {\it Transonic shock in a nozzle I: 2D case.} Comm. Pure Appl. Math.,58 (8) (2005), 999-1050.

\end{thebibliography}

\end{document}